\documentclass[12pt]{article}
\usepackage{amsmath,amsfonts,amssymb,amsthm}
\input amssym.def
\begin{document}
\def \Z{\Bbb Z}
\def \C{\Bbb C}
\def \R{\Bbb R}
\def \Q{\Bbb Q}
\def \N{\Bbb N}

\def \A{{\mathcal{A}}}
\def \D{{\mathcal{D}}}
\def \E{{\mathcal{E}}}
\def \E{{\mathcal{E}}}
\def \H{\mathcal{H}}
\def \S{{\mathcal{S}}}
\def \V{{\mathcal{V}}ir}
\def \wt{{\rm wt}}
\def \tr{{\rm tr}}
\def \Gr{{\rm Gr}}
\def \Cl{{\rm Cl}}
\def \Map{{\rm Map}}
\def \span{{\rm span}}
\def \Res{{\rm Res}}
\def \Der{{\rm Der}}
\def \End{{\rm End}}
\def \Ind {{\rm Ind}}
\def \Irr {{\rm Irr}}
\def \Aut{{\rm Aut}}
\def \GL{{\rm GL}}
\def \Hom{{\rm Hom}}
\def \mod{{\rm mod}}
\def \ann{{\rm Ann}}
\def \ad{{\rm ad}}
\def \rank{{\rm rank}\;}
\def \<{\langle}
\def \>{\rangle}

\def \g{{\frak{g}}}
\def \h{{\hbar}}
\def \k{{\frak{k}}}
\def \sl{{\frak{sl}}}
\def \gl{{\frak{gl}}}

\def \be{\begin{equation}\label}
\def \ee{\end{equation}}
\def \bex{\begin{example}\label}
\def \eex{\end{example}}
\def \bl{\begin{lem}\label}
\def \el{\end{lem}}
\def \bt{\begin{thm}\label}
\def \et{\end{thm}}
\def \bp{\begin{prop}\label}
\def \ep{\end{prop}}
\def \br{\begin{rem}\label}
\def \er{\end{rem}}
\def \bc{\begin{coro}\label}
\def \ec{\end{coro}}
\def \bd{\begin{de}\label}
\def \ed{\end{de}}

\newcommand{\m}{\bf m}
\newcommand{\n}{\bf n}
\newcommand{\nno}{\nonumber}
\newcommand{\nord}{\mbox{\scriptsize ${\circ\atop\circ}$}}
\newtheorem{thm}{Theorem}[section]
\newtheorem{prop}[thm]{Proposition}
\newtheorem{coro}[thm]{Corollary}
\newtheorem{conj}[thm]{Conjecture}
\newtheorem{example}[thm]{Example}
\newtheorem{lem}[thm]{Lemma}
\newtheorem{rem}[thm]{Remark}
\newtheorem{de}[thm]{Definition}
\newtheorem{hy}[thm]{Hypothesis}
\makeatletter \@addtoreset{equation}{section}
\def\theequation{\thesection.\arabic{equation}}
\makeatother \makeatletter

\begin{center}
{\Large \bf  $G$-covariant $\phi$-coordinated quasi modules for quantum vertex algebras}
\end{center}

\begin{center}
{Haisheng Li\footnote{Partially supported by NSA grant
H98230-11-1-0161 and China NSF grant (No. 11128103)}\\
Department of Mathematical Sciences\\
Rutgers University, Camden, NJ 08102, USA}
\end{center}

\begin{abstract}
This is a paper in a series to study quantum vertex algebras and their relations with various quantum algebras.
In this paper, we introduce a notion of T-type quantum vertex algebra and a notion of $G$-covariant $\phi$-coordinated
quasi module for a $T$-type quantum vertex algebra with an automorphism group $G$. We refine and extend several previous results
and we obtain a commutator formula for $G$-covariant $\phi$-coordinated quasi modules.
As an illustrating example, we study a special case of the deformed Virasoro algebra $\V_{p,q}$ with $q=-1$, to which
we associate a Clifford vertex superalgebra and its $G$-covariant $\phi$-coordinated quasi modules.
\end{abstract}

\section{Introduction}

Previously, inspired by Etingof-Kazhdan's theory of quantum vertex operator algebras (see \cite{ek}),
we developed a theory of (weak) quantum vertex algebras and their modules in \cite{li-qva1} and \cite{li-qva2}.
It was shown  (see \cite{li-qva2}, \cite{kl}) that quantum vertex algebras in this sense can be associated to
Zamolodchikov-Faddeev algebras of a certain type.
To associate more (quantum) algebras including quantum affine algebras with quantum vertex algebras,
we furthermore developed in \cite{li-phi-module} a theory of what we called $\phi$-coordinated quasi modules
for a weak quantum vertex algebra. Indeed, this new theory enables us to associate quantum affine algebras
with weak quantum vertex algebras {\em conceptually}. It then remains to realize this conceptual association
{\em explicitly} and to study the associated quantum vertex algebras and their modules, including their $\phi$-coordinated quasi modules.

In this paper, we develop the theory of $\phi$-coordinated quasi modules further and
we refine and extend some of the results therein, paving the way to explicitly associate
quantum vertex algebras to quantum affine algebras.
Specifically, we introduce a notion of $T$-type quantum vertex algebra and a notion of
$G$-covariant $\phi$-coordinated quasi module for a $T$-type quantum vertex algebra $V$,
where $G$ is an automorphism group of $V$ equipped with a linear character $\chi$.
In \cite{li-phi-module}, a Jacobi-type identity for $\phi$-coordinated modules over a weak quantum vertex algebra was obtained.
Among the main results of this paper, we establish a commutator formula for $G$-covariant $\phi$-coordinated {\em quasi} modules
for a T-type weak quantum vertex algebra.
As an illustrating example, we apply this theory to the deformed Virasoro algebra $\V_{p,q}$ with $q=-1$
and we obtain a natural connection of this algebra with a concrete Clifford vertex superalgebra.

Now, we give a more detailed account of this paper.
First, the notion of quantum vertex algebra formulated in \cite{li-qva1} was somewhat derived from
Etingof-Kazhdan's notion of quantum vertex operator algebra. Quantum vertex algebras in this sense are
generalizations of vertex algebras and vertex super-algebras, whereas
quantum vertex operator algebras in the sense of Etingof-Kazhdan
are formal deformations of vertex algebras.
Among the main ingredients for a quantum vertex algebra $V$ is a
rational quantum Yang-Baxter operator $\S(x)$ which governs a generalized commutativity (namely locality)
of the associated vertex operators on $V$. On the one hand,
this notion singles out of a large family of vertex algebra-like structures, and on the other hand,
this theory is not much more complicated than the ordinary vertex algebra theory.

In regard to the theory of $\phi$-coordinated quasi modules
for weak quantum vertex algebras, this very parameter $\phi$ is what was called an associate of
the one-dimensional additive formal group (law) $F(x,y)=x+y$,
where an associate is a formal series $\phi(x,z)\in \C((x))[[z]]$, satisfying
\begin{eqnarray}
\phi(x,0)=x \ \mbox{ and }\ \phi(x,\phi(y,z))=\phi(x+y,z)\ \ (=\phi(F(x,y),z)).
\end{eqnarray}
It was shown that every associate (of the additive formal group law) can be explicitly constructed by
$\phi(x,z)=e^{zp(x)\frac{d}{dx}}\cdot x$ with $p(x)\in \C((x))$. In particular, taking $p(x)=1$ one gets
$\phi(x,z)=x+z=F(x,z)$ (the formal group itself) and taking $p(x)=x$ one gets
$\phi(x,z)=xe^{z}$.
The truth is that the ordinary theory of (weak) quantum vertex algebras and modules is governed by the formal group law, whereas
to each associate $\phi$ one can attach a theory of $\phi$-coordinated quasi modules.
In the theory of weak quantum vertex algebras and modules, an important role was played by the associativity
\begin{eqnarray}
Y(u,x_{0}+x_{2})Y(v,x_{2})=Y(Y(u,x_{0})v,x_{2})
\end{eqnarray}
({\em unrigorous}), which is commonly referred to as operator product expansion by physicists.
In contrast, in the theory of $\phi$-coordinated (quasi) modules, this associativity is replaced with
\begin{eqnarray}\label{ephi-associativity}
Y(u,\phi(x_{2},x_{0}))Y(v,x_{2})=Y(Y(u,x_{0})v,x_{2})
\end{eqnarray}
({\em unrigorous}).
For most of the better known algebras including quantum affine algebras, we take $\phi(x,z)=xe^{z}$.
It was proved therein that for any highest weight module $W$
for a quantum affine algebra of a fixed level, the Drinfeld generating functions,
viewed as fields on $W$, generate a weak quantum vertex algebra in a certain natural way
with $W$ as a canonical $\phi$-coordinated quasi module. In this way,
quantum affine algebras were tied to weak quantum vertex algebras conceptually.

Recall that a rational quantum Yang-Baxter operator is an important ingredient for a quantum vertex algebra.
In this paper, we define a {\em T-type quantum vertex algebra} to be a quantum vertex algebra $V$
such that the rational quantum Yang-Baxter operator $\S(x)$ is associated to a trigonometric quantum Yang-Baxter operator $T(x)$
with $\S(x)=T(e^{x})$.
Let $V$ be a T-type quantum vertex algebra and let $G$ be an automorphism group equipped with
a linear character $\chi: G\rightarrow \C^{\times}$. We then define a
notion of $G$-covariant $\phi$-coordinated quasi $V$-module.
In the definition, in addition to the associativity (\ref{ephi-associativity})
with $\phi(x_{2},x_{0})=x_{2}e^{x_{0}}$, two defining properties are that
$$Y_{W}(gv,x)=Y_{W}(v,\chi(g)x)\ \ \mbox{ for }g\in G, \ v\in V$$
and that for any $u,v\in V$, there exists a polynomial $p(x)$ with only zeroes in $\{ \chi(g)\ |\ g\in G\}$
such that
$$p(x/z)Y_{W}(u,x)Y_{W}(v,z)\in \Hom (W,W((x,z))).$$

There is a similarity between this notion and the notion of twisted $V$-module (see \cite{flm}, \cite{ffr}, \cite{dong})
for a vertex algebra $V$, where a twisted $V$-module
is associated to an automorphism $\sigma$ of finite order $N$ with
a choice of an $N$-th root of unity, which amounts to choosing a linear character for $\<\sigma\>=\Z_{N}$.
(For a vertex operator algebra $V$, a notion of $\Z_{N}$-covariant quasi module was introduced  and
a canonical isomorphism between the category of $\Z_{N}$-twisted modules and the category of
$\Z_{N}$-covariant quasi modules was established in \cite{li-tqe}.)
As one of the main results of this paper, we obtain a commutator formula, which is somewhat analogous to that for twisted modules.

In this paper, as an illustrating example we apply this theory to the deformed Virasoro algebra $\V_{p,q}$ with $q=-1$.
The deformed Virasoro algebra $\V_{p,q}$ with nonzero complex parameters $p$ and $q$ was introduced in \cite{skao} and
the special case with $q=-1$ was later studied in \cite{bp}.
The deformed Virasoro algebra is among the simplest quantum algebras
in the sense that it is generated by one single field.
We here establish a canonical isomorphism between the category of highest weight modules for
$\V_{p,q}$ with $q=-1$ and the category of $G$-covariant $\phi$-coordinated quasi modules
for a specific Clifford vertex super-algebra.

As for quantum affine algebras $U_{q}(\hat{\g})$ (with $q$ a nonzero complex number), our speculation
is that for every $\ell\in \C$, we shall have a canonical T-type quantum vertex algebra
with an automorphism group $G$ isomorphic to a subgroup of $\C^{\times}$
such that highest weight modules for $U_{q}(\hat{\g})$
of level $\ell$ are exactly $G$-covariant $\phi$-coordinated quasi modules for the conjectured $T$-type quantum vertex algebra.

In a pioneer work \cite{efr},  E. Frenkel and N. Reshetikhin introduced a theory of deformed chiral algebras.
Among the key ingredients in the notion of deformed chiral algebra are two spaces, called the space of states and the space of fields, and
a {\em trigonometric} quantum Yang-Baxter operator on the space of fields.
This theory and the Etingof-Kazhdan theory look quite different and they have been studied independently.
Now, the structure of a deformed chiral algebra can be enhanced to a $T$-type quantum vertex algebra $V$
equipped with an automorphism group $G$ together with a $G$-covariant $\phi$-coordinated quasi module.

This paper is organized as follows: Sections 2 and 3 are preliminary; In Section 2 we present some basic results on formal calculus, and 
 in Section 3 we recall the basic notions and results about
(weak) quantum vertex algebras. In Section 4, we define the notions of $T$-type quantum vertex algebra and $G$-covariant $\phi$-coordinated quasi modules, and we present several results including the commutator formula. In Section 5, we study
the deformed Virasoro algebra $\V_{p,q}$ with $q=-1$.

\section{Some results on formal calculus}
In this section, we shall review some formal variable notations and conventions, and
we formulate (and prove) certain basic results which we need in later sections.

First of all, throughout this paper all vector spaces are assumed to be over $\C$ (the field of complex numbers) and
$\N$ denotes the set of nonnegative integers. Letters $x,y,z, x_{i},y_{i},z_{i}$ with $i=0,1,\dots$
will be mutually commuting independent formal variables.

For a vector space $U$,  $U[[x,x^{-1}]]$ denotes the space of formal infinite integer power series
with coefficients in $U$,
$U[[x]]$ denotes the space of formal infinite nonnegative power series, and
$U((x))$ denotes the space of lower truncated formal Laurent series.
In particular, we have a vector space $\C[[x,x^{-1}]]$, an algebra $\C[[x]]$, and a field $\C((x))$.

We denote by $F_{\C[[x]]}$ the field of fractions of $\C[[x]]$, which is naturally isomorphic to $\C((x))$.
Similarly, denote by $F_{\C[[x_{1},x_{2}]]}$ the field of fractions of $\C[[x_{1},x_{2}]]$. We define
\begin{eqnarray}
\iota_{x_{1},x_{2}}: F_{\C[[x_{1},x_{2}]]}\rightarrow \C((x_{1}))((x_{2}))
\end{eqnarray}
to be the canonical extension of the embedding of $\C[[x_{1},x_{2}]]$ into $\C((x_{1}))((x_{2}))$.
Following \cite{flm} and \cite{fhl}, we use the following formal variable convention:
For $n\in \Z$,
\begin{eqnarray}
(x_{1}-x_{2})^{n}=\iota_{x_{1},x_{2}}\left((x_{1}-x_{2})^{n}\right)=\sum_{i\ge 0}\binom{n}{i}(-1)^{i}x_{1}^{n-i}x_{2}^{i}.
\end{eqnarray}
For $A(x)=\sum_{n\in \Z}A(n)x^{n}\in U[[x,x^{-1}]]$ with $U$ a vector space, we have
\begin{eqnarray*}
&&e^{z\frac{d}{dx}}A(x)=A(x+z)=\sum_{n\in \Z}A(n)(x+z)^{n},\\
&&e^{z(x\frac{d}{dx})}A(x)=A(xe^{z}).
\end{eqnarray*}

We also use the following formal delta functions:
$$\delta(x)=\sum_{n\in \Z}x^{n},$$
$$x_{1}^{-1}\delta\left(\frac{x_{2}}{x_{1}}\right)=\sum_{n\in \Z}x_{1}^{-n-1}x_{2}^{n},$$
$$x_{0}^{-1}\delta\left(\frac{x_{1}-x_{2}}{x_{0}}\right)=\sum_{n\in \Z}x_{0}^{-n-1}(x_{1}-x_{2})^{n}
=\sum_{n\in \Z}\sum_{i\ge 0}\binom{n}{i}(-1)^{i}x_{0}^{-n-1}x_{1}^{n-i}x_{2}^{i}.$$
For $n\ge 0$, we have
\begin{eqnarray}
\frac{1}{n!}\left(\frac{\partial}{\partial x_{2}}\right)^{n}x_{1}^{-1}\delta\left(\frac{x_{2}}{x_{1}}\right)
=(x_{1}-x_{2})^{-n-1}-(-x_{2}+x_{1})^{-n-1}.
\end{eqnarray}
An important property is the following substitution rule:
\begin{eqnarray}
f(x_{1},x_{2})\delta\left(\frac{x_{2}}{x_{1}}\right)=f(x_{2},x_{2})\delta\left(\frac{x_{2}}{x_{1}}\right),
\end{eqnarray}
where $f(x_{1},x_{2})$ is any series such that $f(x_{2},x_{2})$ exists, e.g., $f(x_{1},x_{2})\in U((x_{1},x_{2}))$.

We shall frequently use the following simple facts:

\bl{lannihilating-property} Let $\lambda\in \C^{\times}$ and let $k$ be a positive integer. Then
\begin{eqnarray}
(x_{1}-\lambda x_{2})^{k}\left(x_{2}\frac{\partial}{\partial x_{2}}\right)^{j}\delta\left(\frac{\lambda x_{2}}{x_{1}}\right)=0
\end{eqnarray}
for $0\le j\le k-1$.
\el

\begin{proof} Using another independent formal variable $z$, we have
\begin{eqnarray*}
&&(x_{1}-\lambda x_{2})^{k}\sum_{j\ge 0}\frac{z^{j}}{j!}\left(x_{2}\frac{\partial}{\partial x_{2}}\right)^{j}\delta\left(\frac{\lambda x_{2}}{x_{1}}\right)\nonumber\\
&=&(x_{1}-\lambda x_{2})^{k}e^{zx_{2}\frac{\partial}{\partial x_{2}}}\delta\left(\frac{\lambda x_{2}}{x_{1}}\right)\nonumber\\
&=&(x_{1}-\lambda x_{2})^{k}\delta\left(\frac{\lambda x_{2}e^{z}}{x_{1}}\right)\nonumber\\
&=&(\lambda x_{2})^{k}(e^{z}-1)^{k}\delta\left(\frac{\lambda x_{2}e^{z}}{x_{1}}\right)\nonumber\\
&\in &z^{k}\C[[x_{1}^{\pm 1},x_{2}^{\pm 1},z]].
\end{eqnarray*}
 Then it follows immediately.
\end{proof}

\bl{lexpression-unique} Let $\lambda_{1},\dots,\lambda_{r}$ be distinct nonzero complex numbers,
let $k$ be a nonnegative integer, and let
$A_{ij}(x)\in U[[x,x^{-1}]]$ for $1\le i\le r,\ 0\le j\le k$, where $U$ is a vector space. Then
\begin{eqnarray}\label{esum=0}
\sum_{i=1}^{r}\sum_{j=0}^{k} A_{ij}(x_{2})\left(x_{2}\frac{\partial}{\partial x_{2}}\right)^{j}
\delta\left(\frac{\lambda_{i}x_{2}}{x_{1}}\right)=0
\end{eqnarray}
if and only if $A_{ij}(x)=0$ for all $i,j$.
\el

\begin{proof} We only need to prove the ``only if'' part.
For $1\le i\le r$, set
$$p_{i}(x)=\frac{(x-\lambda_{1})^{k+1}\cdots (x-\lambda_{r})^{k+1}}{(x-\lambda_{i})^{k+1}}\in \C[x].$$
Then there exist $q_{1}(x),\dots,q_{r}(x)\in \C[x]$ such that
\begin{eqnarray}\label{eexpression-of1}
p_{1}(x)q_{1}(x)+\cdots +p_{r}(x)q_{r}(x)=1.
\end{eqnarray}
Let $1\le s\le r$ be arbitrarily fixed.
For $1\le i\le r$ with $i\ne s$, as $(x-\lambda_{i})^{k+1}$ divides $p_{s}(x)$,
by Lemma \ref{lannihilating-property} we have
$$p_{s}(x_{1}/x_{2})\left(x_{2}\frac{\partial}{\partial x_{2}}\right)^{j}
\delta\left(\frac{\lambda_{i}x_{2}}{x_{1}}\right)=0 \ \ \ \mbox{ for }0\le j\le k.$$
Then multiplying both sides of (\ref{esum=0}) by $p_{s}(x_{1}/x_{2})q_{s}(x_{1}/x_{2})$, we get
\begin{eqnarray*}
\sum_{j=0}^{k} A_{sj}(x_{2})p_{s}(x_{1}/x_{2})q_{s}(x_{1}/x_{2})\left(x_{2}\frac{\partial}{\partial x_{2}}\right)^{j}
\delta\left(\frac{\lambda_{s}x_{2}}{x_{1}}\right)=0.
\end{eqnarray*}
Noticing that for $1\le n\le r$ with $n\ne s$ and for $0\le j\le k$,
$$p_{n}(x_{1}/x_{2})q_{n}(x_{1}/x_{2})\left(x_{2}\frac{\partial}{\partial x_{2}}\right)^{j}
\delta\left(\frac{\lambda_{s}x_{2}}{x_{1}}\right)=0,$$
 then using (\ref{eexpression-of1}) we obtain
\begin{eqnarray*}
\sum_{j=0}^{k} A_{sj}(x_{2})\left(x_{2}\frac{\partial}{\partial x_{2}}\right)^{j}
\delta\left(\frac{\lambda_{s}x_{2}}{x_{1}}\right)=0,
\end{eqnarray*}
which amounts to
\begin{eqnarray}
\sum_{j=0}^{k} A_{sj}(x_{2})\left(x_{2}\frac{\partial}{\partial x_{2}}\right)^{j}
\delta\left(\frac{x_{2}}{x_{1}}\right)=0.
\end{eqnarray}
For any nonzero integer $m$, by extracting the coefficient of $x_{1}^{-m}$ we get
\begin{eqnarray*}
\sum_{j=0}^{k} m^{j}A_{sj}(x_{2})=0.
\end{eqnarray*}
Taking $m=1,2,\dots,k+1$ and solving the system of equations we obtain
$A_{sj}(x)=0$ for $0\le j\le k$, as desired.
\end{proof}

Let $W$ be a vector space. Set
\begin{eqnarray}
\E(W)=\Hom (W,W((x)))\subset (\End W)[[x,x^{-1}]].
\end{eqnarray}

\bl{lformal-locality} Let $a(x),b(x)\in \E(W),\ K(x_{1},x_{2})\in \Hom (W,W((x_{2}))((x_{1})))$,
$$p(x)=(x-\lambda_{1})^{k_{1}}\cdots (x-\lambda_{r})^{k_{r}}\in \C[x],$$
where $\lambda_{1},\dots,\lambda_{r}$ are distinct nonzero complex numbers and
$k_{1},\dots,k_{r}$ are positive integers. Then
\begin{eqnarray}
p(x_{1}/x_{2})\left(a(x_{1})b(x_{2})-K(x_{1},x_{2})\right)=0
\end{eqnarray}
if and only if
\begin{eqnarray}
a(x_{1})b(x_{2})-K(x_{1},x_{2})=\sum_{i=1}^{r}\sum_{j=0}^{k_{i}-1}A_{ij}(x_{2}) \left(x_{2}\frac{\partial}{\partial x_{2}}\right)^{j}\delta\left(\lambda_{i}\frac{x_{2}}{x_{1}}\right)
\end{eqnarray}
for some $A_{ij}(x)\in \E(W)$, which are uniquely determined.
\el

\begin{proof} The uniqueness follows from Lemma \ref{lexpression-unique},
so it remains to prove the existence. {}From assumption we have
$$p(x_{1}/x_{2})a(x_{1})b(x_{2})=p(x_{1}/x_{2})K(x_{1},x_{2}),$$
which implies
$$p(x_{1}/x_{2})a(x_{1})b(x_{2}),\ \ p(x_{1}/x_{2})K(x_{1},x_{2})\in \Hom (W,W((x_{1},x_{2}))).$$
Thus
\begin{eqnarray*}
&&a(x_{1})b(x_{2})=\iota_{x_{1},x_{2}}\left(\frac{1}{p(x_{1}/x_{2})}\right)A(x_{1},x_{2}),\\
& &K(x_{1},x_{2})=\iota_{x_{2},x_{1}}\left(\frac{1}{p(x_{1}/x_{2})}\right)A(x_{1},x_{2})
\end{eqnarray*}
for some $A(x_{1},x_{2})\in \Hom(W,W((x_{1},x_{2})))$. Then
\begin{eqnarray*}
a(x_{1})b(x_{2})-K(x_{1},x_{2})=\left(\iota_{x_{1},x_{2}}\left(\frac{1}{p(x_{1}/x_{2})}\right)
-\iota_{x_{2},x_{1}}\left(\frac{1}{p(x_{1}/x_{2})}\right)\right)A(x_{1},x_{2}).
\end{eqnarray*}
Write
\begin{eqnarray*}
\frac{1}{p(x)}=\sum_{i=1}^{r}\sum_{j=0}^{k_{i}}\frac{a_{ij}}{(x-\lambda_{i})^{j}}
\end{eqnarray*}
with $a_{ij}\in \C$. We then have
\begin{eqnarray}
&&a(x_{1})b(x_{2})-K(x_{1},x_{2})\nonumber\\
&=&A(x_{1},x_{2})\left(\iota_{x_{1},x_{2}}\left(\frac{1}{p(x_{1}/x_{2})}\right)
-\iota_{x_{2},x_{1}}\left(\frac{1}{p(x_{1}/x_{2})}\right)\right)\nonumber\\
&=&A(x_{1},x_{2})\sum_{i=1}^{r}\sum_{j=1}^{k_{i}}\left(\frac{a_{ij}x_{2}^{j}}{(x_{1}-\lambda_{i}x_{2})^{j}}
-\frac{a_{ij}x_{2}^{j}}{(-\lambda_{i}x_{2}+x_{1})^{j}}\right)\nonumber\\
&=&\sum_{i=1}^{r}\sum_{j=1}^{k_{i}}a_{ij}\frac{1}{(j-1)!}\lambda_{i}^{1-j}x_{2}^{j}A(x_{1},x_{2})\left(\frac{\partial}{\partial x_{2}}\right)^{j-1}x_{1}^{-1}\delta\left(\frac{\lambda_{i}x_{2}}{x_{1}}\right).
\end{eqnarray}
It follows from induction that for any nonnegative integer $n$ and for any
$B(x_{1},x_{2})\in \Hom (W,W((x_{1},x_{2})))$, we have
\begin{eqnarray}
B(x_{1},x_{2})\left(\frac{\partial}{\partial x_{2}}\right)^{n}x_{1}^{-1}\delta\left(\frac{\lambda_{i}x_{2}}{x_{1}}\right)
=\sum_{i=0}^{n}B_{i}(x_{2})\left(\frac{\partial}{\partial x_{2}}\right)^{i}\delta\left(\frac{\lambda_{i}x_{2}}{x_{1}}\right)
\end{eqnarray}
for some $B_{i}(x)\in \E(W)$. On the other hand, for $n\ge 1$ we have
\begin{eqnarray}
\left(\frac{\partial}{\partial x_{2}}\right)^{n}=
x_{2}^{-n}\left(x_{2}\frac{\partial}{\partial x_{2}}\right)\left(x_{2}\frac{\partial}{\partial x_{2}} -1\right)\cdots \left(x_{2}\frac{\partial}{\partial x_{2}}-(n-1)\right).
\end{eqnarray}
Using all these we get
\begin{eqnarray*}
a(x_{1})b(x_{2})-K(x_{1},x_{2})=\sum_{i=1}^{r}\sum_{j=0}^{k_{i}-1}A_{ij}(x_{2})\left(x_{2}\frac{\partial}{\partial x_{2}}\right)^{j}
\delta\left(\frac{\lambda_{i} x_{2}}{x_{1}}\right)
\end{eqnarray*}
for some $A_{ij}(x)\in \E(W)$, as desired.
\end{proof}

We shall also need the following simple fact:

\bl{lsimple-formal-cal} Let $W$ be a vector space, let $A(x_{1},x_{2})\in \Hom (W,W((x_{1},x_{2})))$, and let $\lambda\in \C^{\times}$.
If $A(x_{1},x_{2})\ne 0$, then
$$A(x_{1},x_{2})=(x_{1}-\lambda x_{2})^{k}B(x_{1},x_{2})$$
for some $k\in \N,\ B(x_{1},x_{2})\in \Hom (W,W((x_{1},x_{2})))$ with $B(\lambda x_{2},x_{2})\ne 0$.
\el

\begin{proof} Notice that for any $m\in \Z$,
$$x_{1}^{m}-(\lambda x_{2})^{m}=(x_{1}-\lambda x_{2})F_{m}(x_{1},x_{2})$$
for some $F_{m}(x_{1},x_{2})\in \C[x_{1}^{\pm 1},x_{2}^{\pm 1}]$.
For any $C(x_{1},x_{2})\in \Hom (W,W((x_{1},x_{2})))$, if $C(\lambda x_{2},x_{2})=0$, then
$$C(x_{1},x_{2})=C(x_{1},x_{2})-C(\lambda x_{2},x_{2})=(x_{1}-\lambda x_{2})\bar{C}(x_{1},x_{2})$$
for some $\bar{C}(x_{1},x_{2})\in \Hom (W,W((x_{1},x_{2})))$. On the other hand, as $A(x_{1},x_{2})\ne 0$,
$A(x_{1},x_{2})w\ne 0$ for some $w\in W$, where $A(x_{1},x_{2})w\in W((x_{1},x_{2}))$. Then
$A(x_{1},x_{2})w=x_{1}^{r}x_{2}^{s}G(x_{1},x_{2})$ for some $r,s\in \Z,\ G(x_{1},x_{2})\in W[[x_{1},x_{2}]]$ (nonzero).
Assuming $G(x_{1}+\lambda x_{2},x_{2})=\sum_{m\ge 0}G_{m}(x_{2})x_{1}^{m}$ with $G_{m}(x_{2})\in \C[[x_{2}]]$, we have
$$G(x_{1},x_{2})=\sum_{m\ge 0}G_{m}(x_{2})(x_{1}-\lambda x_{2})^{m}.$$
It follows that
$G(x_{1},x_{2})=(x_{1}-\lambda x_{2})^{n}H(x_{1},x_{2})$ with $n\in \N$
and $H(x_{1},x_{2})\in W[[x_{1},x_{2}]]$ such that $H(\lambda x_{2},x_{2})\ne 0$.
Then
$$A(x_{1},x_{2})w=(x_{1}-\lambda x_{2})^{n}x_{1}^{r}x_{2}^{s}H(x_{1},x_{2}).$$
In view of these, there exists a largest nonnegative integer $k$ such that
$A(x_{1},x_{2})=(x_{1}-\lambda x_{2})^{k}B(x_{1},x_{2})$ for some
$B(x_{1},x_{2})\in \Hom (W,W((x_{1},x_{2})))$ with $B(\lambda x_{2},x_{2})\ne 0$.
\end{proof}

\section{Quantum vertex algebras and their modules}
In this section, we recall the notions of
(weak) quantum vertex algebras and their modules, including a conceptual construction.

We first recall the notion of weak quantum vertex algebra,
which was formulated and studied in \cite{li-qva1} and \cite{li-qva2}.

\bd{dweak-qva} {\em A {\em weak quantum vertex algebra} is a vector
space $V$ equipped with a linear map
\begin{eqnarray*}
Y(\cdot,x):&& V\rightarrow \Hom(V,V((x)))\subset (\End
V)[[x,x^{-1}]]\\
&&v\mapsto Y(v,x)=\sum_{n\in \Z}v_{n}x^{-n-1}\ \ (\mbox{where
}v_{n}\in \End V),
\end{eqnarray*}
called the {\em adjoint vertex operator map}, and a vector ${\bf
1}\in V$, called the {\em vacuum vector}, satisfying the following
conditions: For $v\in V$,
$$Y({\bf 1},x)v=v,\ \ Y(v,x){\bf 1}\in V[[x]]\ \ \mbox{and }\
\lim_{x\rightarrow 0}Y(v,x){\bf 1}=v, $$
and for $u,v\in V$, there exist
$$u^{(i)},v^{(i)}\in V,\ \ f_{i}(x)\in \C((x))\ \ \mbox{ for
}i=1,\dots,r$$ such that
\begin{eqnarray}\label{sjacobi-identity}
&&x_{0}^{-1}\delta\left(\frac{x_{1}-x_{2}}{x_{0}}\right)Y(u,x_{1})Y(v,x_{2})\nonumber\\
&&\hspace{1cm}-x_{0}^{-1}\delta\left(\frac{x_{2}-x_{1}}{-x_{0}}\right)
\sum_{i=1}^{r}f_{i}(x_{2}-x_{1})Y(v^{(i)},x_{2})Y(u^{(i)},x_{1})\nonumber\\
&=&x_{2}^{-1}\delta\left(\frac{x_{1}-x_{0}}{x_{2}}\right)Y(Y(u,x_{0})v,x_{2})
\end{eqnarray}
(the {\em $\S$-Jacobi identity}).} \ed

The following was proved in \cite{li-qva1}:

\bp{pwqva-equivalence} In Definition \ref{dweak-qva},
the $\S$-Jacobi identity axiom can be equivalently replaced by {\em weak associativity}:
 For $u,v,w\in V$, there exists a nonnegative integer $l$ such that
\begin{eqnarray}
(x_{0}+x_{2})^{l}Y(u,x_{0}+x_{2})Y(v,x_{2})w=(x_{0}+x_{2})^{l}Y(Y(u,x_{0})v,x_{2})w,
\end{eqnarray}
 and {\em $\S$-locality:} For any $u,v\in V$, there exist
 $$u^{(i)},v^{(i)}\in V,\ \ f_{i}(x)\in \C((x))\ \ \mbox{ for
}i=1,\dots,r,$$ and a nonnegative integer $k$ such that
\begin{eqnarray}\label{es-locality-relation}
(x_{1}-x_{2})^{k}Y(u,x_{1})Y(v,x_{2})=(x_{1}-x_{2})^{k}\sum_{i=1}^{r}f_{i}(x_{2}-x_{1})Y(v^{(i)},x_{2})Y(u^{(i)},x_{1}).
\end{eqnarray}
\ep

Let $V$ be a weak quantum vertex algebra.
Define a linear operator $\D$ on $V$ by $\D (v)=v_{-2}{\bf 1}$ for $v\in V$. Then
\begin{eqnarray}
[\D, Y(v,x)]=Y(\D (v),x)=\frac{d}{dx}Y(v,x)
\end{eqnarray}
for $v\in V$. It was proved in \cite{li-qva1} that
the $\S$-locality relation (\ref{es-locality-relation}) amounts to
\begin{eqnarray}
Y(u,x)v=\sum_{i=1}^{r}f_{i}(-x)e^{x\D}Y(v^{(i)},-x)u^{(i)}
\end{eqnarray}
(the {\em $\S$-skew symmetry}).

\bd{dmodule-wqva}
{\em Let $V$ be a weak quantum vertex algebra. A {\em $V$-module} is a vector space $W$ equipped with a linear map
$$Y_{W}(\cdot,x): \ V\rightarrow \Hom (W,W((x)))\subset (\End W)[[x,x^{-1}]],$$
satisfying the conditions that $Y_{W}({\bf 1},x)=1_{W}$ (the identity operator on $W$) and that for any $u,v\in V,\ w\in W$,
there exists a nonnegative integer $l$ such that
\begin{eqnarray}
(x_{0}+x_{2})^{l}Y_{W}(u,x_{0}+x_{2})Y_{W}(v,x_{2})w=(x_{0}+x_{2})^{l}Y_{W}(Y(u,x_{0})v,x_{2})w.
\end{eqnarray}}
\ed

We have (see \cite{li-qva1}):

\bp{pwqva-module-con}
Let $V$ be a weak quantum vertex algebra and let $(W,Y_{W})$ be a $V$-module. Then
\begin{eqnarray}
Y_{W}(\D v,x)=\frac{d}{dx}Y_{W}(v,x)\ \ \ \mbox{ for }v\in V.
\end{eqnarray}
For $u,v\in V$, let $u^{(i)},v^{(i)}\in V, \ f_{i}(x)\in \C((x))$ $(i=1,\dots,r)$
such that $\S$-Jacobi identity (\ref{sjacobi-identity}) holds. Then
\begin{eqnarray}
&&x_{0}^{-1}\delta\left(\frac{x_{1}-x_{2}}{x_{0}}\right)Y_{W}(u,x_{1})Y_{W}(v,x_{2})\nonumber\\
&&\hspace{1cm}-x_{0}^{-1}\delta\left(\frac{x_{2}-x_{1}}{-x_{0}}\right)
\sum_{i=1}^{r}f_{i}(x_{2}-x_{1})Y_{W}(v^{(i)},x_{2})Y_{W}(u^{(i)},x_{1})\nonumber\\
&=&x_{2}^{-1}\delta\left(\frac{x_{1}-x_{0}}{x_{2}}\right)Y_{W}(Y(u,x_{0})v,x_{2}).
\end{eqnarray}
\ep

We next recall from \cite{li-qva1} the conceptual construction of weak quantum vertex algebras and modules.

\bd{dcompatibility}
{\em Let $W$ be a vector space and let $a(x),b(x)\in \E(W)$. The ordered pair $(a(x),b(x))$ is said to be {\em compatible}
if there exists $k\in \N$ such that
\begin{eqnarray}\label{ecompatibility}
(x_{1}-x_{2})^{k}a(x_{1})b(x_{2})\in \Hom (W,W((x_{1},x_{2}))).
\end{eqnarray}}
\ed

Let $(a(x),b(x))$ be a compatible pair in $\E(W)$.
Define $a(x)_{n}b(x)\in \E(W)$ for $n\in \Z$ in terms of generating function
$$Y_{\E}(a(x),z)b(x)=\sum_{n\in \Z}a(x)_{n}b(x) z^{-n-1}$$
by
\begin{eqnarray}
Y_{\E}(a(x),z)b(x)=z^{-k}\left( (x_{1}-x)^{k}a(x_{1})b(x)\right)|_{x_{1}=x+z},
\end{eqnarray}
where $k$ is any nonnegative integer such that (\ref{ecompatibility}) holds.

\bd{ds-locality}
{\em A subset $U$ of $\E(W)$ is said to be {\em $\S$-local} if for any $a(x),b(x)\in U$, there exist
$$a^{(i)}(x),b^{(i)}(x)\in U,\ f_{i}(x)\in \C((x))\ \ (i=1,\dots,r)$$ such that
\begin{eqnarray}
(x_{1}-x_{2})^{k}a(x_{1})b(x_{2})=(x_{1}-x_{2})^{k}\sum_{i=1}^{r}f_{i}(x_{2}-x_{1})b^{(i)}(x_{2})a^{(i)}(x_{1})
\end{eqnarray}
for some nonnegative integer $k$.}
\ed

Note that if $U$ is an $\S$-local subset of $\E(W)$, then every pair $(a(x),b(x))$ in $U$ is compatible.
The following result was obtained in \cite{li-qva1}:

\bt{tmain-wqva-construction}
Every $\S$-local subset $U$ of $\E(W)$ generates a weak quantum vertex algebra $\<U\>$
with $W$ as a canonical module with $Y_{W}(\alpha(x),z)=\alpha(z)$ for $\alpha(x)\in \<U\>$.
\et

Let $U$ be a vector space. Recall that a {\em rational quantum Yang-Baxter operator} on $U$ is a linear map
$$\S(x): U\otimes U\rightarrow U\otimes U\otimes \C((x)),$$
satisfying
\begin{eqnarray}
\S_{12}(x)\S_{13}(x+z)\S_{23}(z)=\S_{23}(z)\S_{13}(x+z)\S_{12}(x)
\end{eqnarray}
(the {\em quantum Yang-Baxter equation}), where for $1\le i<j\le 3$,
$$\S^{ij}(x):U\otimes U\otimes U\rightarrow U\otimes U\otimes
U\otimes \C((x))$$ denotes the canonical extension of $\S(x)$.  It
is said to be {\em unitary} if
$$\S(x)\S^{21}(-x)=1,$$
where $\S^{21}(x)=\sigma \S(x)\sigma$ with $\sigma$ denoting the
flip operator on $U\otimes U$.

For a weak quantum vertex algebra $V$, following \cite{ek} let
\begin{eqnarray}
Y(x): V\otimes V\rightarrow V((x))
\end{eqnarray}
be the canonical linear map associated to the vertex operator map
$Y(\cdot,x)$.

The following notion, which was formulated in \cite{li-qva1},
was derived from Etingof-Kazhdan's notion of quantum vertex operator algebra (see \cite{ek}):

\bd{dqva} {\em A {\em quantum vertex algebra} is a weak quantum
vertex algebra $V$ equipped with a unitary rational quantum
Yang-Baxter operator $\S(x)$ on $V$, satisfying the following
conditions:
\begin{eqnarray}
&&\S(x)({\bf 1}\otimes v)={\bf 1}\otimes v\ \ \ \mbox{ for }v\in
V,\label{esvacuum}\\
&&[\D\otimes 1, \S(x)]=-\frac{d}{dx}\S(x),\label{d1s}\\
&&Y(u,x)v=e^{x\D}Y(-x)\S(-x)(v\otimes u)\ \ \mbox{ for }u,v\in V,\\
&&\S(x_{1})(Y(x_{2})\otimes 1)=(Y(x_{2})\otimes
1)\S^{23}(x_{1})\S^{13}(x_{1}+x_{2}),\label{sy1}
\end{eqnarray}
where $\D$ is the linear operator on $V$ defined by $\D (v)=v_{-2}{\bf 1}$ for $v\in V$.
We sometimes denote a quantum vertex algebra by a pair $(V,\S)$.} \ed

The following notion was due to Etingof-Kazhdan (see \cite{ek}):

\bd{dnon-degeneracy}
{\em Let $V$ be a weak quantum vertex algebra. For any positive integer $n$, define a linear map
\begin{eqnarray*}
Z_{n}: V^{\otimes n}\otimes \C((x_{1}))((x_{2}))\cdots ((x_{n}))\rightarrow V((x_{1}))((x_{2}))\cdots ((x_{n}))
\end{eqnarray*}
by
$$Z_{n}(v^{(1)}\otimes \cdots \otimes v^{(n)}\otimes f)=f Y(v^{(1)},x_{1})\cdots Y(v^{(n)},x_{n}){\bf 1}$$
for $v^{(1)},\dots,v^{(n)}\in V,\ f\in \C((x_{1}))((x_{2}))\cdots ((x_{n}))$.
$V$ is said to be {\em non-degenerate} if for every positive integer $n$, $Z_{n}$ is injective.}
\ed

The following result can be found in \cite{li-qva1} (cf. \cite{ek}):

\bp{pnon-degenerate} Let $V$ be a weak quantum vertex algebra.
Assume that $V$ is non-degenerate. Then there exists a linear map
$\S(x): V\otimes V\rightarrow V\otimes V\otimes \C((x))$, which is
uniquely determined by
\begin{eqnarray*}
Y(u,x)v=e^{x\D}Y(-x)\S(-x)(v\otimes u) \ \ \ \mbox{for }u,v\in V.
\end{eqnarray*}
Furthermore, $(V,\S)$ carries the structure of a quantum vertex
algebra and the following relation holds
\begin{eqnarray}
[1\otimes \D, \S(x)]=\frac{d}{dx}\S(x).
\end{eqnarray}
\ep

\br{rnotions} {\em Note that a quantum vertex algebra was defined as
a pair $(V,\S)$. In view of Proposition \ref{pnon-degenerate}, the
term ``a non-degenerate quantum vertex algebra'' without reference
to a quantum Yang-Baxter operator is unambiguous. If a weak quantum vertex algebra
$V$ is of countable dimension over $\C$ and if $V$ as a $V$-module
is irreducible, then $V$ is non-degenerate by a result of \cite{li-qva2} (Corollary 3.10).
In view of this, the term ``irreducible quantum
vertex algebra'' is also unambiguous.} \er

We conclude this section with some basic notions.
The notions of homomorphism, isomorphism and automorphism for weak quantum vertex algebras
are defined in the obvious way.
For example, an {\em automorphism} of a weak quantum vertex algebra $V$ is
a bijective linear endomorphism $\sigma$ of
$V$ such that $\sigma({\bf 1})={\bf 1}$ and
$$\sigma (Y(u,x)v)=Y(\sigma(u),x)\sigma(v)\ \ \mbox{ for }u,v\in V.$$

\bd{dgenerating-subset} {\em A subset $U$ of a weak quantum vertex algebra $V$ 
is called a {\em generating subset} if $V$ is linearly spanned
by vectors
$$u^{(1)}_{n_{1}}\cdots u^{(r)}_{n_{r}}{\bf 1}$$
for $r\in \N,\ u^{(i)}\in U,\ n_{i}\in \Z$ $(1\le i\le r)$.}
\ed

\section{$\phi$-coordinated quasi modules for weak quantum vertex algebras}

In this section, we first recall from \cite{li-phi-module} the definition of a $\phi$-coordinated quasi module
for a weak quantum vertex algebra and the conceptual construction.
We then define the notions of $T$-type weak quantum vertex algebra and $G$-covariant $\phi$-coordinated
quasi module. As one of the main results we establish a commutator formula.

We begin with a convention. Note that for any $p(x)\in \C[x]$, $p(e^{x})\in \C[[x]]$. One can show
 (cf. \cite{li-phi-module}) that $p(e^{x})=0$ if and only if $p(x)=0$.
Recall that $\C(x)$ denotes the field of rational functions.
For any $g(x)=p(x)/q(x)\in \C(x)$ with $p(x),q(x)\in \C[x]$, we define
\begin{eqnarray}
g(e^{x})=\frac{p(e^{x})}{q(e^{x})}\in F_{\C[[x]]}=\C((x)).
\end{eqnarray}

\bd{dphi-coordinated-quasi-module} {\em Let $V$ be a weak quantum vertex algebra.
A {\em $\phi$-coordinated quasi $V$-module} is a vector space $W$ equipped with a linear map
$$Y_{W}(\cdot,x):\ V\rightarrow \Hom (W,W((x)))\subset (\End W)[[x,x^{-1}]],$$
satisfying the conditions that $Y_{W}({\bf 1},x)=1_{W}$ and that for $u,v\in V$,
there exists a nonzero polynomial $p(x)$ such that
\begin{eqnarray}
p(x_{1}/x_{2})Y_{W}(u,x_{1})Y_{W}(v,x_{2})\in \Hom(W,W((x_{1},x_{2})))
\end{eqnarray}
and
\begin{eqnarray}\label{epphi-associativity}
p(e^{z})Y_{W}(Y(u,z)v,x_{2})=\left(p(x_{1}/x_{2})Y_{W}(u,x_{1})Y_{W}(v,x_{2})\right)|_{x_{1}=x_{2}e^{z}}.
\end{eqnarray}}
\ed

\br{rabout-phi} {\em The parameter $\phi$ in Definition \ref{dphi-coordinated-quasi-module}
 refers to the formal series $\phi(x,z)=xe^{z}$, which is
a particular associate of
the one-dimensional additive formal group (law) $F(x,y)=x+y$, as defined in \cite{li-phi-module}. }
\er

The following was proved in \cite{li-phi-module} (Proposition 5.6):

\bp{precall} Let $V$ be a weak quantum vertex algebra, let $u,v\in V$, and
let $(W,Y_{W})$ be a $\phi$-coordinated quasi $V$-module.
Suppose that
\begin{eqnarray}\label{ekuv-trig}
&&(x_{1}-x_{2})^{k}Y(u,x_{1})Y(v,x_{2})\nonumber\\
&=& (x_{1}-x_{2})^{k}
\sum_{i=1}^{r}\iota_{x_{2},x_{1}}(f_{i}(e^{x_{1}-x_{2}}))Y(v^{(i)},x_{2})Y(u^{(i)},x_{1}),
\end{eqnarray}
where $k\in \N,\ f_{i}(x)\in \C(x),\ u^{(i)},v^{(i)}\in V$,
 and suppose that $p(x)$ is any nonzero polynomial such that
$$p(x_{1}/x_{2})Y_{W}(u,x_{1})Y_{W}(v,x_{2})\in \Hom (W,W((x_{1},x_{2}))).$$
Then
\begin{eqnarray}
&&p(x_{1}/x_{2})Y_{W}(u,x_{1})Y_{W}(v,x_{2})\nonumber\\
&=&p(x_{1}/x_{2})\sum_{i=1}^{r}\iota_{x_{2},x_{1}}\left(f_{i}(x_{1}/x_{2})\right)Y_{W}(v^{(i)},x_{2})Y_{W}(u^{(i)},x_{1}).
\end{eqnarray}
\ep

Set
\begin{eqnarray}
\log (1+z)=\sum_{n\ge 1}(-1)^{n-1}\frac{z^{n}}{n}\in z\C[[z]].
\end{eqnarray}
Notice that $\sum_{n\ge 1}(-1)^{n-1}\frac{z^{n-1}}{n}$ is a unit in $\C[[z]]$.
For $m\in \Z$, it is understood that
$$\left(\log (1+z)\right)^{m}=z^{m}\left(\sum_{n\ge 1}(-1)^{n-1}\frac{z^{n-1}}{n}\right)^{m}\in z^{m}\C[[z]]\subset \C((z)).$$

The following was proved in \cite{li-phi-module} (Lemma 5.8)\footnote{There is a typo in \cite{li-phi-module}: $B(x_{2},x_{1})$ in (5.18) and in the proof should be $B(x_{1},x_{2})$.}:

\bl{lthree-delta-log}  Let $W$ be a vector space and let
$$A(x_{1},x_{2})\in \Hom (W,W((x_{1}))((x_{2}))),\ \ B(x_{1},x_{2})\in \Hom (W,W((x_{2}))((x_{1}))),$$
$$C(x_{0},x_{2})\in \Hom (W,W((x_{2}))((x_{0}))).$$
If there exists a nonnegative integer $k$ such that
\begin{eqnarray*}
&&(x_{1}-x_{2})^{k}A(x_{1},x_{2})=(x_{1}-x_{2})^{k}B(x_{1},x_{2}),\\
&&\left((x_{1}-x_{2})^{k}A(x_{1},x_{2})\right)|_{x_{1}=x_{2}e^{x_{0}}}=x_{2}^{k}(e^{x_{0}}-1)^{k}C(x_{0},x_{2}),
\end{eqnarray*}
then
\begin{eqnarray}
&&(zx_{2})^{-1}\delta\left(\frac{x_{1}-x_{2}}{zx_{2}}\right)A(x_{1},x_{2})
-(zx_{2})^{-1}\delta\left(\frac{x_{2}-x_{1}}{-zx_{2}}\right)B(x_{1},x_{2})\nonumber\\
&&\hspace{1cm}=x_{1}^{-1}\delta\left(\frac{x_{2}(1+z)}{x_{1}}\right)C(\log (1+z),x_{2}).
\end{eqnarray}
Furthermore, the converse is also true.
\el

The following result generalizes Proposition 5.9 of \cite{li-phi-module}:

\bp{pquasiphi-jacobi-identity}
Let $V$ be a weak quantum vertex algebra
and let $(W,Y_{W})$ be a $\phi$-coordinated quasi $V$-module. Let $u,v\in V$. Assume
$$u^{(i)},v^{(i)}\in V,\ f_{i}(x)\in \C(x)\ \ (i=1,\dots,r)$$
such that (\ref{ekuv-trig}) holds for some nonnegative integer $k$.
Then there exists a nonzero polynomial $p(x)$  such that
\begin{eqnarray}
&&(xz)^{-1}\delta\left(\frac{x_{1}-x}{xz}\right)p(x_{1}/x)Y_{W}(u,x_{1})Y_{W}(v,x)\hspace{2cm}\nonumber\\
&&\ \ -(xz)^{-1}\delta\left(\frac{x-x_{1}}{-xz}\right)p(x_{1}/x)\sum_{i=1}^{r}\iota_{x_{1},x}
\left(f_{i}(x_{1}/x)\right)Y_{W}(v^{(i)},x)Y_{W}(u^{(i)},x_{1})\nonumber\\
&=&x_{1}^{-1}\delta\left(\frac{x(1+z)}{x_{1}}\right)p(x_{1}/x)Y_{W}\left(Y(u,\log (1+z))v,x\right).
\end{eqnarray}
\ep

\begin{proof} {}From definition, there exists a nonzero polynomial $p(x)$  such that
\begin{eqnarray}\label{eproof-pywuv}
p(x_{1}/x_{2})Y_{W}(u,x_{1})Y_{W}(v,x_{2})\in \Hom(W,W((x_{1},x_{2})))
\end{eqnarray}
and
\begin{eqnarray}
p(e^{z})Y_{W}(Y(u,z)v,x_{2})=\left(p(x_{1}/x_{2})Y_{W}(u,x_{1})Y_{W}(v,x_{2})\right)|_{x_{1}=x_{2}e^{z}}.
\end{eqnarray}
With (\ref{eproof-pywuv}), by Proposition \ref{precall} we have
\begin{eqnarray}
&&p(x_{1}/x_{2})Y_{W}(u,x_{1})Y_{W}(v,x_{2})\nonumber\\
&=&p(x_{1}/x_{2})\sum_{i=1}^{r}\iota_{x_{2},x_{1}}\left(f_{i}(x_{1}/x_{2})\right)Y_{W}(v^{(i)},x_{2})Y_{W}(u^{(i)},x_{1}).
\end{eqnarray}
Then it follows immediately from Lemma \ref{lthree-delta-log} (with $k=0$).
\end{proof}

Motivated by Proposition \ref{precall} we formulate the following notion:

\bd{dwqva-trig} {\em A weak quantum vertex algebra $V$ is said to be of {\em T-type} if for any $u,v\in V$, there exist
$$u^{(i)},v^{(i)}\in V,\ f_{i}(x)\in \C(x)\ \ (i=1,\dots,r)$$ such that (\ref{ekuv-trig}) holds for some nonnegative integer $k$. }
\ed

Recall that a {\em trigonometric  quantum Yang-Baxter operator} on a vector space $U$ is a linear map
$$\S(x): U\otimes U\rightarrow U\otimes U\otimes \C(x),$$
satisfying quantum Yang-Baxter equation
\begin{eqnarray}
\S_{12}(x)\S_{13}(xz)\S_{23}(z)=\S_{23}(z)\S_{13}(xz)\S_{12}(x).
\end{eqnarray}
It is said to be {\em unitary} if $\S_{21}(1/x)\S(x)=1$.

\bd{dqva-trig} {\em A quantum vertex algebra $(V,\S(x))$ is said to be of {\em T-type} if there exists
a unitary trigonometric quantum Yang-Baxter operator $T(x)$ on $V$ such that $\S(x)=T(e^{x})$. }
\ed

The following technical result follows from the
proof of Proposition 2.6 in \cite{ltw} with $\C((x))$ replaced by $\C(x)$ in a few places:

\bl{lspecial-rational}
Let $V$ be a weak quantum vertex algebra. Suppose that $U$ is a generating subset of $V$ such that
the condition in Definition \ref{dwqva-trig} with $U$ in place of $V$ holds.
 Then $V$ is of T-type.
\el

Next, we recall from \cite{li-phi-module} the conceptual construction of weak quantum vertex algebras and their
$\phi$-coordinated quasi modules.
Let $W$ be a vector space which is fixed for a while. Let $a(x),b(x)\in \E(W)$. Assume that there
exists a nonzero polynomial $p(x)$ such that
\begin{eqnarray}\label{equasi-compatibility}
p(x/z)a(x)b(z)\in \Hom (W,W((x,z))).
\end{eqnarray}
Define $a(x)_{n}^{e}b(x)\in \E(W)$ for $n\in \Z$ in terms of
generating function $$Y_{\E}^{e}(a(x),z)b(x)=\sum_{n\in
\Z}(a(x)_{n}^{e}b(x))z^{-n-1}$$ by
\begin{eqnarray}
Y_{\E}^{e}(a(x),z)b(x)=p(e^{z})^{-1}
\left(p(x_{1}/x)a(x_{1})b(x)\right)|_{x_{1}=xe^{z}},
\end{eqnarray}
where $p(x)$ is any nonzero polynomial such that (\ref{equasi-compatibility}) holds
and $p(e^{z})^{-1}$ denotes the inverse of $p(e^{z})$ in
$\C((z))$.

\bd{dquasi-trig-local} {\em
A subset $U$ of $\E(W)$ is said to be {\em quasi $\S_{trig}$-local} if for any
$a(x),b(x)\in U$, there exist
$$a^{(i)}(x), b^{(i)}(x)\in U,\; f_{i}(x)\in \C(x)\ \ (i=1,\dots,r)$$
such that
\begin{eqnarray}\label{e-puc-comm}
p(x/z)a(x)b(z)=p(x/z)\sum_{i=1}^{r}\iota_{z,x}(f(z/x))b^{(i)}(z)a^{(i)}(x)
\end{eqnarray}
for some nonzero polynomial $p(x)$. A subset $U$ of $\E(W)$ is {\em
quasi $\S_{trig}$-compatible} if for any $a(x),b(x)\in U$, there
exists a nonzero polynomial $p(x)$ such that
\begin{eqnarray*}
p(x/z)a(x)b(z)\in \Hom (W,W((x,z))).
\end{eqnarray*}}
\ed

Note that relation (\ref{e-puc-comm}) implies  (\ref{equasi-compatibility}), so that
every quasi $\S_{trig}$-local subset is automatically quasi $\S_{trig}$-compatible.

The following strengthens one of the main results in \cite{li-phi-module}:

\bt{tmain-phi}
Let $W$ be a vector space and let $U$ be any quasi $\S_{trig}$-local
subset of $\E(W)$. Then $U$ generates a  weak quantum vertex algebra $\<U\>_{e}$ with $W$
as a faithful $\phi$-coordinated quasi module
where $Y_{W}(\alpha(x),z)=\alpha(z)$ for $\alpha(x)\in \<U\>_{e}$.
Furthermore, $\<U\>_{e}$ is of T-type.
\et

\begin{proof} The first part was proved in \cite{li-phi-module} (Theorem 5.4).
Then we have a weak quantum vertex algebra $\<U\>_{e}$ with $U$ as a generating subset.
Let $a(x),b(x)\in U$. There exist
$$a^{(i)}(x), b^{(i)}(x)\in U,\; f_{i}(x)\in \C(x)\ \ (i=1,\dots,r)$$
such that
\begin{eqnarray*}
p(x/z)a(x)b(z)=p(x/z)\sum_{i=1}^{r}\iota_{z,x}(f(z/x))b^{(i)}(z)a^{(i)}(x)
\end{eqnarray*}
for some nonzero polynomial $p(x)$. By Proposition 5.3 of \cite{li-phi-module}, we have
\begin{eqnarray*}
&&(x_{1}-x_{2})^{k}Y_{\E}^{e}(a(x),x_{1})Y_{\E}^{e}(b(x),x_{2})\\
&=&(x_{1}-x_{2})^{k}\sum_{i=1}^{r}\iota_{x_{2},x_{1}}(f_{i}(e^{x_{2}-x_{1}}))Y_{\E}^{e}(b^{(i)}(x),x_{2})Y_{\E}^{e}(a^{(i)}(x),x_{1})
\end{eqnarray*}
for some nonnegative integer $k$.
Then it follows from Lemma \ref{lspecial-rational} that $\<U\>_{e}$ is of T-type.
\end{proof}

We have the following result (cf. \cite{li-phi-module}, Proposition 4.11):

\bp{pquasi-representation}
Let $V$ be a T-type weak quantum vertex algebra and let $(W,Y_{W})$ be a $\phi$-coordinated quasi
$V$-module. Set $V_{W}=\{ Y_{W}(v,x)\ |\ v\in V\}$. Then $V_{W}$ is a quasi $\S_{trig}$-local subspace of $\E(W)$,
$(V_{W},Y_{\E}^{e},1_{W})$ is a weak quantum vertex algebra, and  $Y_{W}$ is a homomorphism of weak quantum vertex algebras.
\ep

\begin{proof} Combining Definition \ref{dphi-coordinated-quasi-module} with Proposition \ref{precall}, we see
that $V_{W}$ is a quasi $\S_{trig}$-local subspace of $\E(W)$. Let $u,v\in V$.
By definition, there exists a nonzero polynomial $p(x)$ such that
$$p(x_{1}/x)Y_{W}(u,x_{1})Y_{W}(v,x)\in \Hom(W,W((x_{1},x)))$$
and
$$p(e^{z})Y_{W}(Y(u,z)v,x)=\left(p(x_{1}/x)Y_{W}(u,x_{1})Y_{W}(v,x)\right)|_{x_{1}=xe^{z}}.$$
On the other hand, from the definition of $Y_{\E}^{e}(\cdot,x)$ we have
$$p(e^{z})Y_{\E}^{e}(Y_{W}(u,x),z)Y_{W}(v,x)=\left(p(x_{1}/x)Y_{W}(u,x_{1})Y_{W}(v,x)\right)|_{x_{1}=xe^{z}}.$$
Consequently,
$$p(e^{z})Y_{\E}^{e}(Y_{W}(u,x),z)Y_{W}(v,x)=p(e^{z})Y_{W}(Y(u,z)v,x).$$
Since both $Y_{\E}^{e}(Y_{W}(u,x),z)Y_{W}(v,x)$ and $Y_{W}(Y(u,z)v,x)$ involve only finitely many negative powers of $z$,
by cancellation we obtain
$$Y_{\E}^{e}(Y_{W}(u,x),z)Y_{W}(v,x)=Y_{W}(Y(u,z)v,x).$$
It follows that $(V_{W},Y_{\E}^{e},1_{W})$ is a weak quantum vertex algebra and
$Y_{W}$ is a homomorphism of weak quantum vertex algebras.
\end{proof}

Note that $a(\lambda x)\in \E(W)$ for any $\lambda\in \C^{\times},\ a(x)\in \E(W)$.
For $\lambda\in \C^{\times}$, define $R_{\lambda}\in \End (\E(W))$ by
\begin{eqnarray}
R_{\lambda}(a(x))=a(\lambda x)\ \ \ \mbox{ for }a(x)\in \E(W).
\end{eqnarray}
This gives rise
to a group action of $\C^{\times}$ on $\E(W)$. {}From now on, we fix this particular action.

\bp{pautomorphism}
Let $W$ be a vector space, let $\Gamma$ be a subgroup of $\C^{\times}$, and let $U$ be a
quasi $\S_{trig}$-local subset of $\E(W)$ such that $U$ is $\Gamma$-stable.
Then $\Gamma$ is an automorphism group of the weak quantum vertex algebra $\<U\>_{e}$ generated by $U$.
\ep

\begin{proof} Let $a(x),b(x)\in \<U\>_{e}$. Then there exists a nonzero polynomial $p(x)$ such that
$$p(x_{1}/x)a(x_{1})b(x)\in \Hom (W,W((x_{1},x))).$$
For any nonzero complex number $\lambda$ we have
$$p(x_{1}/x)a(\lambda x_{1})b(\lambda x)\in \Hom (W,W((x_{1},x))).$$
Then we get
\begin{eqnarray*}
&&Y_{\E}^{e}(a(\lambda x),z)b(\lambda x)\\
&=&p(e^{z})^{-1}\left(p(x_{1}/x)a(\lambda x_{1})b(\lambda x)\right)|_{x_{1}=xe^{z}}\\
&=&\left[p(e^{z})^{-1}\left(p(\bar{x}_{1}/\bar{x})a(\bar{x}_{1})b(\bar{x})\right)|_{\bar{x}_{1}=\bar{x}e^{z}}\right]|_{\bar{x}=\lambda x}\\
&=&\left[Y_{\E}^{e}(a(\bar{x}),z)b(\bar{x})\right]_{\bar{x}=\lambda x}.
\end{eqnarray*}
It follows from induction that $\<U\>_{e}$ is $\Gamma$-stable and $\Gamma$ acts on weak quantum vertex algebra
$\<U\>_{e}$ by automorphisms.
\end{proof}

\bl{lextension}
Let $W$ be a vector space, let $a(x),b(x)\in \E(W)$, and let
$$K(x_{1},x_{2})\in \Hom (W,W((x_{2}))((x_{1}))).$$ Assume
\begin{eqnarray}
a(x_{1})b(x_{2})-K(x_{1},x_{2})&=&\sum_{j=0}^{r}A_{j}(x_{2})\frac{1}{j!}\left(x_{2}\frac{\partial}{\partial x_{2}}\right)^{j}\delta\left(\frac{x_{2}}{x_{1}}\right)\nonumber\\
&&+\sum_{i=1}^{k}\sum_{j=0}^{s}B_{ij}(x_{2})\frac{1}{j!}\left(x_{2}\frac{\partial}{\partial x_{2}}\right)^{j}\delta\left(\frac{\lambda_{i}x_{2}}{x_{1}}\right),
\end{eqnarray}
where $A_{j}(x), B_{ij}(x)\in \E(W)$, $\lambda_{i}\in \C^{\times}$ with $\lambda_{i}\ne 1$ for $1\le i\le k$. Then
\begin{eqnarray}
&&a(x)_{j}^{e}b(x)=A_{j}(x)\ \ \ \mbox{ for }0\le j\le r,\nonumber\\
&&a(x)_{j}^{e}b(x)=0\ \ \ \mbox{ for }j>r.
\end{eqnarray}
\el

\begin{proof}  Set
$$q(x)=\prod_{i=1}^{k}\left((x-1)^{r+1}-(\lambda_{i}-1)^{r+1}\right)^{s+1}\in \C[x].$$
Noticing that
$$q(x)=\bar{q}(x)\prod_{i=1}^{k}(x-\lambda_{i})^{s+1}$$
for some $\bar{q}(x)\in \C[x]$, we have
$$q(x_{1}/x_{2})\left(x_{2}\frac{\partial}{\partial x_{2}}\right)^{j}\delta\left(\frac{\lambda_{i}x_{2}}{x_{1}}\right)=0$$
for $0\le j\le s,\ 1\le i\le k$. Then we get
\begin{eqnarray}\label{eqab=}
q(x_{1}/x_{2})\left(a(x_{1})b(x_{2})-K(x_{1},x_{2})\right)
=\sum_{j=0}^{r}q(x_{1}/x_{2})A_{j}(x_{2})\frac{1}{j!}\left(x_{2}\frac{\partial}{\partial x_{2}}\right)^{j}\delta\left(\frac{x_{2}}{x_{1}}\right).
\end{eqnarray}
This implies
$$(x_{1}/x_{2}-1)^{r+1}q(x_{1}/x_{2})(a(x_{1})b(x_{2})-K(x_{1},x_{2}))=0.$$
{}From definition we have
\begin{eqnarray*}
(e^{x_{0}}-1)^{r+1}q(e^{x_{0}})Y_{\E}^{e}(a(x),x_{0})b(x)=\left((x_{1}/x-1)^{r+1}q(x_{1}/x)a(x_{1})b(x)\right)|_{x_{1}=xe^{x_{0}}}.
\end{eqnarray*}
With the above two identities, by Lemma \ref{lthree-delta-log} we obtain
\begin{eqnarray*}
&&x_{1}^{-1}\delta\left(\frac{x(1+z)}{x_{1}}\right)q(x_{1}/x)Y_{\E}^{e}(a(x),\log (1+z))b(x)\nonumber\\
&=&(xz)^{-1}\delta\left(\frac{x_{1}-x}{xz}\right)q(x_{1}/x)a(x_{1})b(x)\hspace{2cm}\nonumber\\
&&-(xz)^{-1}\delta\left(\frac{x-x_{1}}{-xz}\right)q(x_{1}/x)K(x_{1},x).
\end{eqnarray*}
Applying $\Res_{z}$ to both sides, we get
\begin{eqnarray}
&&q(x_{1}/x)(a(x_{1})b(x)-K(x_{1},x))\nonumber\\
&=&\Res_{x_{0}} q(x_{1}/x)Y_{\E}^{e}(a(x),x_{0})b(x)e^{x_{0}x\frac{\partial}{\partial x}}\delta\left(\frac{x}{x_{1}}\right)\nonumber\\
&=&\sum_{i\ge 0}a(x)_{i}^{e}b(x) \frac{1}{i!}q(x_{1}/x)\left(x\frac{\partial}{\partial x}\right)^{i}\delta\left(\frac{x}{x_{1}}\right),
\end{eqnarray}
which is a finite sum.
Combining this with (\ref{eqab=}) we obtain
\begin{eqnarray}\label{efinalcomp}
\sum_{i\ge 0}\left(A_{i}(x)-a(x)_{i}b(x)\right)\frac{1}{i!}q(x_{1}/x)\left(x\frac{\partial}{\partial x}\right)^{i}\delta\left(\frac{x}{x_{1}}\right)=0,
\end{eqnarray}
where we set $A_{i}(x)=0$ for $i>r$.
Write
\begin{eqnarray*}
q(x_{1}/x)=(x_{1}/x-1)^{r+1}P(x_{1}/x)+\alpha,
\end{eqnarray*}
where $P(x)\in \C[x]$ and
$$\alpha=(-1)^{k(s+1)}\prod_{i=1}^{k}(\lambda_{i}-1)^{(r+1)(s+1)}.$$
Note that $\alpha\ne 0$ as $\lambda_{i}\ne 1$ for $1\le i\le k$.
Since
$$(x_{1}/x-1)^{r+1}\left(x\frac{\partial}{\partial x}\right)^{j}\delta\left(\frac{x}{x_{1}}\right)=0
\ \ \ \mbox{ for }0\le j\le r,$$
 (\ref{efinalcomp}) reduces to
\begin{eqnarray}\label{elastcomp}
\sum_{i\ge 0}\left(A_{i}(x)-a(x)_{i}b(x)\right)\frac{1}{i!}\alpha
\left(x\frac{\partial}{\partial x}\right)^{i}\delta\left(\frac{x}{x_{1}}\right)=0,
\end{eqnarray}
which implies
$$\alpha \left(A_{i}(x)-a(x)_{i}b(x)\right)=0\ \ \ \mbox{ for }i\ge 0.$$
 Then our assertions follow immediately.
\end{proof}

Using Lemma \ref{lextension} we have the following generalization:

\bp{pextension}
Let $W$ be a vector space, let $a(x),b(x)\in \E(W)$, and let
$$K(x_{1},x_{2})\in \Hom (W,W((x_{2}))((x_{1}))).$$ Assume
\begin{eqnarray}\label{eab=given}
a(x_{1})b(x_{2})-K(x_{1},x_{2})=
\sum_{i=1}^{k}\sum_{j=0}^{r}A_{ij}(x_{2})\frac{1}{j!}\left(x_{2}\frac{\partial}{\partial x_{2}}\right)^{j}\delta\left(\frac{\lambda_{i}x_{2}}{x_{1}}\right),
\end{eqnarray}
where $A_{ij}(x)\in \E(W)$ and $\lambda_{i}\in \C^{\times}$ distinct for $1\le i\le k$. Then
\begin{eqnarray}
&&a(\lambda_{i}x)_{j}^{e}b(x)=A_{ij}(x)\ \ \ \mbox{ for }0\le j\le r,\nonumber\\
&&a(\lambda_{i}x)_{j}^{e}b(x)=0\ \ \ \mbox{ for }j>r.
\end{eqnarray}
\ep

\begin{proof} For $\lambda\in \C^{\times}$, {}from (\ref{eab=given}) we have
\begin{eqnarray*}
a(\lambda x_{1})b(x_{2})-K(\lambda x_{1},x_{2})=
\sum_{i=1}^{k}\sum_{j=0}^{r}A_{ij}(x_{2})\frac{1}{j!}\left(x_{2}\frac{\partial}{\partial x_{2}}\right)^{j}\delta\left(\frac{\lambda_{i}x_{2}}{\lambda x_{1}}\right).
\end{eqnarray*}
Taking $\lambda=\lambda_{i}$ with $1\le i\le k$, using Lemma \ref{lextension} we get
$$a(\lambda_{i}x)_{j}^{e}b(x)=A_{ij}(x)\ \ \ \mbox{ for } 0\le j\le r$$
and $a(\lambda_{i}x)_{j}^{e}b(x)=0$ for $j>r$,
as desired.
\end{proof}

Combining Lemma \ref{lformal-locality} with Proposition \ref{pextension}  we immediately have:

\bc{cabstract-commutator-formula}
Let $W$ be a vector space, let $a(x),b(x)\in \E(W)$, and let
$$K(x_{1},x_{2})\in \Hom (W,W((x_{2}))((x_{1}))).$$
Assume
\begin{eqnarray}
p(x_{1}/x_{2})a(x_{1})b(x_{2})=p(x_{1}/x_{2})K(x_{1},x_{2}),
\end{eqnarray}
where $p(x)=(x-\lambda_{1})^{k_{1}}\cdots (x-\lambda_{r})^{k_{r}}$ with $\lambda_{1},\dots,\lambda_{r}$ distinct nonzero complex numbers
and with $k_{i}\ge 1$. Then
\begin{eqnarray}
a(x_{1})b(x_{2})-K(x_{1},x_{2})=\Res_{x_{0}}\sum_{i=1}^{r} Y_{\E}^{e}(a(\lambda_{i} x),x_{0})b(x)
e^{x_{0}\left(x\frac{\partial}{\partial x}\right)}\delta\left(\frac{\lambda_{i} x}{x_{1}}\right).
\end{eqnarray}
\ec

In the following we shall establish certain technical results.

\bl{lWregularity} Let $a(x),b(x)\in \E(W)$ and let $\lambda\in \C^{\times}$.
Suppose that $k$ is an integer and $p(x)$ is a
polynomial with $p(\lambda)\ne 0$ such that
$$(x_{1}/x_{2}-\lambda)^{k}p(x_{1}/x_{2})a(x_{1})b(x_{2})\in \Hom (W,W((x_{1},x_{2}))).$$
Then
\begin{eqnarray}
a(\lambda x)_{n}^{e}b(x)=0\ \ \ \mbox{ for all }n\ge k.
\end{eqnarray}
 \el

\begin{proof} {}From assumption we have
$$(x_{1}/x_{2}-1)^{k}p(\lambda x_{1}/x_{2})a(\lambda x_{1})b(x_{2})\in \Hom (W,W((x_{1},x_{2}))).$$
Then
\begin{eqnarray*}
&&z^{k}Y_{\E}^{e}(a(\lambda x),z)b(x)\\
&=&z^{k}(e^{z}-1)^{-k}p(\lambda e^{z})^{-1} \left((x_{1}/x-1)^{k}p(\lambda x_{1}/x)a(\lambda x_{1})b(x)\right)|_{x_{1}=xe^{z}},
\end{eqnarray*}
where $p(\lambda e^{z})^{-1}$ denotes the inverse of $p(\lambda e^{z})$ in $\C((z))$.
We  have
$z^{k}(e^{z}-1)^{-k}\in \C[[z]]$ and  $p(\lambda e^{z})^{-1}\in \C[[z]]$ as $p(\lambda)\ne 0$.
We know that
$$\left((x_{1}/x-1)^{k}p(\lambda x_{1}/x)a(\lambda x_{1})b(x)\right)|_{x_{1}=xe^{z}}$$
contains only nonnegative powers of $z$. Then it follows that $z^{k}Y_{\E}^{e}(a(\lambda x),z)b(x)$
involves only nonnegative powers of $z$. Thus,
$a(\lambda x)_{n}^{e}b(x)=0$ for all $n\ge k$.
\end{proof}

\bl{ltruncation}
Let $a(x),b(x)\in \E(W)$. Suppose that $p(x)$ is a nonzero polynomial such that
$$p(x_{1}/x_{2})a(x_{1})b(x_{2})\in \Hom (W,W((x_{1},x_{2}))).$$
Write $p(x)=(x-1)^{s}q(x)$ with $s\in \N,\ q(x)\in \C[x]$ such that $q(1)\ne 0$. Then
for $k\in \Z$,
\begin{eqnarray}\label{ekqabhomw}
(x_{1}/x_{2}-1)^{k}q(x_{1}/x_{2})a(x_{1})b(x_{2})\in \Hom (W,W((x_{1},x_{2})))
\end{eqnarray}
if and only if $a(x)_{j}^{e}b(x)=0$ for $j\ge k$.
\el

\begin{proof} The ``only if'' part follows from Lemma \ref{lWregularity}.
Now, assume that $a(x)_{j}^{e}b(x)=0$ for $j\ge k$. We have
$$a(x_{1})b(x)=\iota_{x_{1},x}(1/p(x_{1}/x))A(x_{1},x)$$
for some $A(x_{1},x)\in \Hom (W,W((x_{1},x)))$. If $A(x_{1},x)=0$, there is nothing to prove.
Assume $A(x_{1},x)\ne 0$. By Lemma \ref{lsimple-formal-cal}
we have $A(x_{1},x)=(x_{1}/x-1)^{r}B(x_{1},x)$ for some $r\in \N,\ B(x_{1},x)\in \Hom (W,W((x_{1},x)))$
with $B(x,x)\ne 0$.  Then
\begin{eqnarray}\label{etemp-need}
(x_{1}/x-1)^{s-r}q(x_{1}/x)a(x_{1})b(x)=B(x_{1},x)\in \Hom (W,W((x_{1},x))),
\end{eqnarray}
so that
\begin{eqnarray*}
Y_{\E}^{e}(a(x),z)b(x)=(e^{z}-1)^{r-s}q(e^{z})^{-1}B(xe^{z},x).
\end{eqnarray*}
We have
$$\lim_{z\rightarrow 0}z^{s-r}Y_{\E}^{e}(a(x),z)b(x)=\lim_{z\rightarrow 0}\left(\frac{e^{z}-1}{z}\right)^{r-s}q(e^{z})^{-1}B(xe^{z},x)=q(1)B(x,x).$$
We get $a(x)_{n}^{e}b(x)=0$ for $n\ge s-r$  and $a(x)_{s-r-1}^{e}b(x)\ne 0$.
As we are given that $a(x)_{m}^{e}b(x)$ for $m\ge k$, we must have $k\ge s-r$. Then it follows from (\ref{etemp-need})
that (\ref{ekqabhomw}) holds.
\end{proof}

\bd{dquasi-phi-module}
{\em Let $V$ be a T-type weak quantum vertex algebra and let $G$ be an automorphism group equipped with a linear character
 $\chi: G\rightarrow \C^{\times}$. A {\em $G$-covariant $\phi$-coordinated quasi $V$-module}
is a $\phi$-coordinated quasi $V$-module $(W,Y_{W})$,
satisfying the conditions that
\begin{eqnarray}\label{eGcovariance-def}
Y_{W}(gv,x)=Y_{W}(v,\chi(g)x)\ \ \ \mbox{ for }g\in G,\ v\in V,
\end{eqnarray}
and that for $u,v\in V$, there exists $p(x)\in \C[x]$ with only zeroes in $\chi(G)$ such that
\begin{eqnarray}\label{epyuyvhom}
p(x_{1}/x_{2})Y_{W}(u,x_{1})Y_{W}(v,x_{2})\in \Hom (W,W((x_{1},x_{2}))).
\end{eqnarray}}
\ed

As the main result of this section we have:

\bt{tcommutator-quasi-phi-module}
Let $V, G,\chi$ be given as in Definition \ref{dquasi-phi-module}
and let $(W,Y_{W})$ be a $G$-covariant $\phi$-coordinated quasi $V$-module. Suppose that $u,v\in V$,
$$u^{(i)},v^{(i)}\in V,\ f_{i}(x)\in \C(x)\ \ (i=1,\dots,r)$$ such that
\begin{eqnarray*}
(x_{1}-x_{2})^{k}Y(u,x_{1})Y(v,x_{2})=(x_{1}-x_{2})^{k}\sum_{i=1}^{r}\iota_{x_{2},x_{1}}f_{i}(e^{x_{1}-x_{2}})Y(v^{(i)},x_{2})Y(u^{(i)},x_{1})
\end{eqnarray*}
for some nonnegative integer $k$. Then there are finitely many $g_{1},\dots,g_{n}\in G$
with $\chi(g_{1}),\dots, \chi(g_{n})$ distinct such that
\begin{eqnarray}
&&Y_{W}(u,x_{1})Y_{W}(v,x)-\sum_{i=1}^{r}\iota_{x_{1},x}
\left(f_{i}(x_{1}/x)\right)Y_{W}(v^{(i)},x)Y_{W}(u^{(i)},x_{1})\nonumber\\
&=&\Res_{x_{0}}\sum_{j=1}^{n}Y_{W}(Y(g_{j}(u),x_{0})v,x)e^{x_{0}\left(x\frac{\partial}{\partial x}\right)}\delta\left(\frac{\chi(g_{j})x}{x_{1}}\right).
\end{eqnarray}
Furthermore, if the linear character $\chi: G\rightarrow \C^{\times}$ is injective, we have
\begin{eqnarray}
&&Y_{W}(u,x_{1})Y_{W}(v,x)-\sum_{i=1}^{r}\iota_{x_{1},x}
\left(f_{i}(x_{1}/x)\right)Y_{W}(v^{(i)},x)Y_{W}(u^{(i)},x_{1})\nonumber\\
&=&\Res_{x_{0}}\sum_{g\in G}Y_{W}(Y(gu,x_{0})v,x)e^{x_{0}\left(x\frac{\partial}{\partial x}\right)}\delta\left(\frac{\chi(g)x}{x_{1}}\right),
\end{eqnarray}
which is a finite sum.
\et

\begin{proof} From definition, there are distinct nonzero complex numbers $\lambda_{1},\dots, \lambda_{n}\in \chi(G)$ such that
$$p(x_{1}/x_{2})Y_{W}(u,x_{1})Y_{W}(v,x_{2})\in \Hom (W,W((x_{1},x_{2}))),$$
where $p(x)=(x-\lambda_{1})^{k_{1}}\cdots (x-\lambda_{n})^{k_{n}}$ with $k_{i}\ge 1$.
By Proposition \ref{precall} we have
\begin{eqnarray*}
&&p(x_{1}/x_{2})Y_{W}(u,x_{1})Y_{W}(v,x_{2})\\
&=&p(x_{1}/x_{2})\sum_{i=1}^{r}\iota_{x_{1},x_{2}}
\left(f_{i}(x_{1}/x_{2})\right)Y_{W}(v^{(i)},x_{2})Y_{W}(u^{(i)},x_{1}).
\end{eqnarray*}
Then by Corollary \ref{cabstract-commutator-formula} we have
\begin{eqnarray}
&&Y_{W}(u,x_{1})Y_{W}(v,x)-\sum_{i=1}^{r}\iota_{x_{1},x}
\left(f_{i}(x_{1}/x)\right)Y_{W}(v^{(i)},x)Y_{W}(u^{(i)},x_{1})\nonumber\\
&=&\Res_{x_{0}}\sum_{j=1}^{n}Y_{\E}^{e}\left(\overline{Y_{W}(u,\lambda_{j}x)},x_{0}\right)\overline{Y_{W}(v,x)}
e^{x_{0}\left(x\frac{\partial}{\partial x}\right)}\delta\left(\frac{\lambda_{j}x}{x_{1}}\right).
\end{eqnarray}
For $1\le j\le n$,  let $g_{j}\in G$ be such that $\lambda_{j}=\chi(g_{j})$. Then
$$Y_{W}(u,\lambda_{j}x)=Y_{W}(u,\chi(g_{j})x)=Y_{W}(g_{j}u,x).$$
Using the fact that $Y_{W}(\cdot,x)$ is a homomorphism (from Proposition \ref{pquasi-representation}), we have
$$Y_{\E}^{e}\left(\overline{Y_{W}(u,\lambda_{j}x)},x_{0}\right)\overline{Y_{W}(v,x)}
=Y_{\E}^{e}\left(\overline{Y_{W}(g_{j}u,x)},x_{0}\right)\overline{Y_{W}(v,x)}
=Y_{W}(Y(g_{j}u,x_{0})v,x).$$
Then we obtain
\begin{eqnarray*}
&&Y_{W}(u,x_{1})Y_{W}(v,x)-\sum_{i=1}^{r}\iota_{x_{1},x}
\left(f_{i}(x_{1}/x)\right)Y_{W}(v^{(i)},x)Y_{W}(u^{(i)},x_{1})\nonumber\\
&=&\Res_{x_{0}}\sum_{j=1}^{n}Y_{W}(Y(g_{j}(u),x_{0})v,x)
e^{x_{0}\left(x\frac{\partial}{\partial x}\right)}\delta\left(\frac{\lambda_{j}x}{x_{1}}\right),
\end{eqnarray*}
proving the first assertion.

Assume that $\chi$ is injective.
For $g\in G$, if $g\ne g_{j}$ for $1\le j\le n$, with $\chi(g)\ne \lambda_{j}$
we have $p(\chi(g))\ne 0$. By Lemma \ref{lWregularity} we get $Y_{W}(u,\chi(g)x)_{m}^{e}Y_{W}(v,x)=0$ for $m\ge 0$.
Then
$$Y_{W}(Y(gu,x_{0})v,x)=Y_{\E}^{e}(Y_{W}(gu,x),x_{0})Y_{W}(v,x)=Y_{\E}^{e}(Y_{W}(u,\chi(g)x),x_{0})Y_{W}(v,x),$$
involving only nonnegative powers of $x_{0}$.
Thus
$$\Res_{x_{0}}Y_{W}(Y(gu,x_{0})v,x)e^{x_{0}\left(x\frac{\partial}{\partial x}\right)}\delta\left(\frac{\chi(g)x}{x_{1}}\right)=0.$$
Then the second assertion follows.
\end{proof}

As an immediate consequence of Theorem \ref{tcommutator-quasi-phi-module} we have:

\bc{ccompatiblity-property} Let $V$ be a T-type weak quantum vertex algebra
 and let $(W,Y_{W})$ be a $G$-covariant $\phi$-coordinated quasi $V$-module.
Let $u,v\in V$.  Then
\begin{eqnarray}
p(x_{1}/x_{2})Y_{W}(u,x_{1})Y_{W}(v,x_{2})\in \Hom
(W,W((x_{1},x_{2}))),
\end{eqnarray}
where $p(x)$ is a product of $(x-\chi(\sigma))^{k_{\sigma}}$ with $\sigma\in G$ and $k_{\sigma}\ge 1$
  such that
$$(\sigma u)_{k_{\sigma}-1}v\ne 0\ \ \mbox{ and }\ \ (\sigma u)_{n}v=0\ \ \mbox{ for } n\ge k_{\sigma}.$$
 \ec

The following is a useful technical result:

\bl{l-phicoordinated-module-induction}
Let $V$ be a $T$-type weak quantum vertex algebra, let $G$ be an automorphism group with a linear character $\chi$,
 and let $(W,Y_{W})$ be a $\phi$-coordinated
quasi $V$-module. Suppose that $U$ is a generating subset of $V$ such that
\begin{eqnarray}\label{ecovariant-lemma}
Y_{W}(gu,x)=Y_{W}(u,\chi(g)x)\ \ \mbox{ for }u\in U,\ g\in G
\end{eqnarray}
and such that $\{ Y_{W}(u,x)\ |\ u\in U\}$ is $\chi(G)$-quasi $\S_{trig}$-local.
Then $(W,Y_{W})$ is a $G$-covariant $\phi$-coordinated quasi $V$-module.
\el

\begin{proof} Set $K=\{ v\in V\ |\ Y_{W}(gv,x)=Y_{W}(v,\chi(g)x)\ \mbox{ for }g\in G\}$. {}From assumption, we have $U\cup \{ {\bf 1}\}\subset K$.
Let $u,v\in K$, $g\in G$. As $Y_{W}(\cdot,x)$ is a homomorphism of weak quantum vertex algebras
by Proposition \ref{pquasi-representation}, we get
\begin{eqnarray*}
Y_{W}(gY(u,z)v,x)&=&Y_{W}(Y(gu,z)gv,x)\\
&=&Y_{\E}^{e}(Y_{W}(gu,x),z)Y_{W}(gv,x)\\
&=&Y_{\E}^{e}(Y_{W}(u,\chi(g)x),z)Y_{W}(v,\chi(g)x)\\
&=&Y_{W}(Y(u,z)v,\chi(g)x)
\end{eqnarray*}
for $g\in G$.
This shows $u_{n}v\in K$ for $n\in \Z$.
Since $U$ generates $V$, we have $V=K$, proving that (\ref{ecovariant-lemma}) holds with $V$ in place of $U$, i.e.,
(\ref{eGcovariance-def}) holds.

Recall Lemma 5.2 from \cite{li-phi-module} that every quasi $S_{trig}$-local subset of
$\E(W)$ is quasi compatible. Let $\Gamma$ be a subgroup of $\C^{\times}$. We say that
a finite sequence $\psi_{1}(x),\dots,\psi_{r}(x)$ in $\E(W)$ is
{\em $\Gamma$-quasi compatible} if there exists a nonzero polynomial
$p(x)$ with only zeros in $\Gamma$ such that
\begin{eqnarray}
\left(\prod_{1\le i<j\le r}p(x_{i}/x_{j})\right)\psi_{1}(x_{1})\cdots \psi_{r}(x_{r})
\in \Hom (W,W((x_{1},\dots,x_{r}))).
\end{eqnarray}
Furthermore, a subset $U$  of $\E(W)$ is
{\em $\Gamma$-quasi compatible} if any finite sequence in $U$ is $\Gamma$-quasi compatible.
The proof of the very Lemma 5.2 with the obvious changes shows that
every $\Gamma$-quasi $S_{trig}$-local subset is $\Gamma$-quasi compatible.

On the other hand, recall Proposition 4.9 in \cite{li-phi-module}: Let $$\psi_{1}(x),\dots,\psi_{r}(x),
a(x),b(x),\phi_{1}(x), \dots,\phi_{s}(x) \in \E(W).$$
Assume that the ordered sequences $(a(x), b(x))$ and
$$(\psi_{1}(x),\dots,\psi_{r}(x), a(x),b(x),\phi_{1}(x),
\dots,\phi_{s}(x))$$
 are quasi compatible.
Then for any $n\in \Z$, the ordered sequence
$$(\psi_{1}(x),\dots,\psi_{r}(x),a(x)_{n}^{e}b(x),\phi_{1}(x),\dots,\phi_{s}(x))$$
is quasi compatible. The same proof with the obvious changes shows that the refinement of Proposition 4.9 with quasi-compatibility replaced by
$\Gamma$-quasi compatibility holds. It then follows from induction that for any $\Gamma$-quasi compatible subset $A$ of $\E(W)$,
$\<A\>_{e}$ is $\Gamma$-quasi compatible.

Note that $Y_{W}(\cdot,x)$ is a homomorphism of weak quantum vertex algebras from $V$ to $V_{W}$
where $V_{W}=\{ Y_{W}(v,x)\ |\ v\in V\}$. Since $U$ is a generating subset of $V$,
it follows that $V_{W}=\<U_{W}\>_{e}$ with $U_{W}=\{ Y_{W}(u,x)\ |\ u\in U\}$.
Then $V_{W}$ is $\chi(G)$-quasi compatible. In particular, (\ref{epyuyvhom}) holds. Therefore,
$(W,Y_{W})$ is a $G$-covariant $\phi$-coordinated quasi $V$-module.
 \end{proof}

As we need, recall the following result from \cite{li-phi-module}
(Lemma 6.7):

\bl{labuv} Let $W$ be a vector space and let $a(x),b(x)\in \E(W)$. Assume that there exist
$$0\ne p(x)\in
\C[x],\; q_{i}(x)\in \C((x)),\;u^{(i)}(x),v^{(i)}(x)\in \E(W) \ \
(1\le i\le r)$$ such that
\begin{eqnarray}\label{epab=quv0}
p(x_{1}/x_{2})a(x_{1})b(x_{2})=\sum_{i=1}^{r}q_{i}(x_{1}/x_{2})u^{(i)}(x_{2})v^{(i)}(x_{1}).
\end{eqnarray}
Then $(a(x),b(x))$ is quasi compatible and
\begin{eqnarray*}
&&p(e^{x_{0}})Y_{\E}^{e}(a(x),x_{0})b(x)\\
&=&\Res_{x_{1}}\left(\frac{1}{x_{1}-xe^{x_{0}}}p(x_{1}/x)a(x_{1})b(x)-
\frac{1}{-xe^{x_{0}}+x_{1}}\sum_{i=1}^{r}q_{i}(x_{1}/x)u^{(i)}(x)v^{(i)}(x_{1})\right).
\end{eqnarray*}
Furthermore, if $k$ is the order of zero of $p(x)$ at $1$, then
$a(x)_{n}^{e}b(x)=0$ for $n\ge k$ and
\begin{eqnarray}\label{ea-1b}
&&\frac{1}{k!}p^{(k)}(1)a(x)_{k-1}^{e}b(x)\nonumber\\
&=&\Res_{x_{1}}\left(\frac{1}{x_{1}-x}p(x_{1}/x)a(x_{1})b(x)-
\frac{1}{-x+x_{1}}q(x_{1}/x)u^{(i)}(x)v^{(i)}(x_{1})\right).\ \
\end{eqnarray}
\el

\section{Deformed Virasoro algebra ${\mathcal{V}}ir_{p,q}$ with $q=-1$ }
In this section, we study the deformed Virasoro algebra ${\mathcal{V}}ir_{p,q}$ with
$q=-1$ in the context of quantum vertex algebras and $\phi$-coordinated quasi modules.
We present a canonical connection of this algebra with a Clifford 
vertex super-algebra and its $G$-covariant $\phi$-coordinated quasi modules.

We first recall {}from \cite{skao} the definition of the deformed Virasoro algebra
$\V_{p,q}$. Let $p$ and $q$ be nonzero complex numbers with $p$ {\em not
a root of $-1$}, so that $p^{n}+1\ne 0$ for every integer $n$. Set
$$t=\frac{q}{p}.$$

\bd{ddvirasoro} {\em The deformed Virasoro algebra $\V_{p,q}$ is
defined to be the associative unital algebra over $\C$ with
generators $T_{n}\; (n\in \Z)$, subject to relations
\begin{eqnarray}\label{edvirasoro-comp}
\sum_{l=0}^{\infty}f_{l}(T_{m-l}T_{n+l}-T_{n-l}T_{m+l})
=-\frac{(1-q)(1-p/q)}{1-p} (p^{m}-p^{-m})\delta_{m+n,0}
\end{eqnarray}
for $m,n\in \Z$, where the coefficients $f_{l}$'s are given by
\begin{eqnarray}
f(z)\equiv \sum_{l= 0}^{\infty}f_{l}z^{l}
=\exp\left(\sum_{n=1}^{\infty}\frac{(1-q^{n})(1-t^{-n})}{1+p^{n}}\frac{z^{n}}{n}\right).
\end{eqnarray}}
\ed

\br{rcompletion} {\em Note that in general, $f_{l}\ne 0$ for
infinitely many $l$, so that the expression on the left-hand side of
(\ref{edvirasoro-comp}) is a genuine infinite sum. Because of this,
the algebra $\V_{p,q}$ is a topological algebra (involving a formal
completion). This can be done rigorously by imitating Frenkel-Zhu's construction of
the universal enveloping algebra $U(V)$ of a vertex operator algebra
$V$ (see \cite{fz}). } \er

Form a generating function
\begin{eqnarray}
T(x)=\sum_{n\in \Z}T_{n}x^{-n}.
\end{eqnarray}
Then the defining relations (\ref{edvirasoro-comp}) can be written
as
\begin{eqnarray}\label{pqvirasoro}
&&f(z/x)T(x)T(z) -f(x/z)T(z)T(x)\nonumber\\
&=&-\frac{(1-q)(1-p/q)}{1-p}
\left[\delta\left(\frac{pz}{x}\right)-\delta\left(\frac{z}{px}\right)\right].
\end{eqnarray}

{}From the defining relations (\ref{edvirasoro-comp}), it can be
readily seen that ${\mathcal{V}}ir_{p,q}$ is a $\Z$-graded algebra
with $\deg T_{n}=-n$ for $n\in \Z$. Then one can study highest
weight ${\mathcal{V}}ir_{p,q}$-modules, in a way similar to that for
the ordinary Virasoro algebra, and in fact, the highest weight
${\mathcal{V}}ir_{p,q}$-modules have been studied in \cite{skao} and
\cite{bp}. One can also define restricted modules, which are more
general than highest weight modules.

The deformed Virasoro algebra ${\mathcal{V}}ir_{p,q}$ with $q=-1$
has been extensively studied in \cite{bp} from an associative algebra point of view.
In this case, one has
$$f(z)=\frac{1+z}{1-z}=1+2\sum_{l=1}^{\infty}z^{l}.$$
The defining relations for ${\mathcal{V}}ir_{p,-1}$ read as
\begin{eqnarray}
[T_{m},T_{n}]+2\sum_{l\ge
0}(T_{m-l}T_{n+l}-T_{n-l}T_{m+l})=-2\left(\frac{1+p}{1-p}\right)
(p^{m}-p^{-m})\delta_{m+n,0}
\end{eqnarray}
for $m,n\in \Z$. In terms of the generating function
$T(x)=\sum_{n\in \Z}T_{n}x^{-n}$, we have
\begin{eqnarray}\label{erelation-generating-q-5}
&&\left(\frac{1+\frac{x_{2}}{x_{1}}}{1-\frac{x_{2}}{x_{1}}}\right)T(x_{1})T(x_{2})
-\left(\frac{1+\frac{x_{1}}{x_{2}}}{1-\frac{x_{1}}{x_{2}}}\right)T(x_{2})T(x_{1})\nonumber\\
&&\ \ \ \hspace{1cm} =-2\left(\frac{1+p}{1-p}\right)
\left(\delta\left(\frac{px_{2}}{x_{1}}\right)
-\delta\left(\frac{p^{-1}x_{2}}{x_{1}}\right)\right).
\end{eqnarray}

We say a ${\mathcal{V}}ir_{p,-1}$-module $W$ is {\em restricted} if for any $w\in W$, $T_{n}w=0$ for $n$ sufficiently large,
namely, $T(x)\in \E(W)$.

Throughout this section, we assume that {\em $p$ is not a root of unity}.
Set
\begin{eqnarray}
\Gamma_{p}=\{ p^{n}\ |\ n\in \Z\}\subset \C^{\times}.
\end{eqnarray}

Next, we introduce a vertex super-algebra.

\bd{dspaceE} {\em Let $E$ be a complex vector space with a basis $\{ e^{(r)}\;|\; r\in
\Z\}$ and equip $E$ with a bilinear form $\<\cdot,\cdot\>$ defined by
\begin{eqnarray}
\<e^{(r)},e^{(s)}\>=2(\delta_{r,s+1}+\delta_{r,s-1}) \ \ \ \mbox{
for }r,s\in \Z.
\end{eqnarray}}
\ed

It can be readily seen that this form is symmetric and
non-degenerate. Set
$$L(E)=E\otimes \C[t,t^{-1}].$$
Let $\ell$ be a complex number. Define a bilinear form
$\<\cdot,\cdot\>_{\ell}$ on $L(E)$ by
\begin{eqnarray}
\<e^{(r)}\otimes t^{m},e^{(s)}\otimes
t^{n}\>_{\ell}=2\ell(\delta_{r,s+1}+\delta_{r,s-1})\delta_{m+n+1,0}
\end{eqnarray}
for $r,s,m,n\in \Z$. This form is symmetric and non-degenerate
whenever $\ell\ne 0$. Denote by $\Cl(L(E),\ell)$ the associated
Clifford algebra, which by definition is the associative unital
algebra generated by vector space $L(E)$, subject to relations
\begin{eqnarray}
(a\otimes t^{m})(b\otimes t^{n})+(b\otimes t^{n})(a\otimes
t^{m})=\<a,b\>\delta_{m+n+1,0}\ell
\end{eqnarray}
for $a,b\in E,\; m,n\in \Z$.

\br{rlie-super} {\em The Clifford algebra  $\Cl(L(E),\ell)$ can
be defined alternatively as the quotient algebra of the universal enveloping
algebra of the Lie super-algebra
$$\hat{E}=(E\otimes \C[t,t^{-1}])\oplus \C {\bf k},$$
modulo relation ${\bf k}=\ell$, where $\hat{E}^{0}=\C {\bf k}$,
$\hat{E}^{1}=E\otimes \C[t,t^{-1}]$, and $[{\bf k}, \hat{E}]=0$,
$$[a\otimes t^{m},b\otimes t^{n}]=\<a,b\>\delta_{m+n+1,0}{\bf k}$$
for $a,b\in E,\; m,n\in \Z$.} \er

Following the tradition, for $a\in E,\; n\in \Z$ we alternatively
denote $a\otimes t^{n}$ by $a_{n}$. For $a\in E$, set
\begin{eqnarray}
a(x)=\sum_{n\in \Z}a_{n}x^{-n-1}\in \Cl(L(E),\ell)[[x,x^{-1}]].
\end{eqnarray}
Set
\begin{eqnarray}
V_{\hat{E}}(\ell,0)= \Cl(L(E),\ell)/\Cl(L(E),\ell)(E\otimes \C[t]),
\end{eqnarray}
a left $\Cl(L(E),\ell)$-module, and set $${\bf 1
}=1+\Cl(L(E),\ell)(E\otimes \C[t])\in V_{\hat{E}}(\ell,0).$$
Identify $E$ as a subspace of $V_{\hat{E}}(\ell,0)$ through the map
$a\mapsto a(-1){\bf 1}$.
It is well known (see \cite{ffr}) that there
exists a vertex super-algebra structure on $V_{\hat{E}}(\ell,0)$,
which is uniquely determined by the condition that ${\bf 1}$ is the
vacuum vector and $Y(a,x)=a(x)$ for $a\in E$.

As we need later, we mention some simple facts about $V_{\hat{E}}(\ell,0)$.
The subspace $E$ generates
$V_{\hat{E}}(\ell,0)$ as a vertex super-algebra. Furthermore, for $r,s\in \Z$, we have
\begin{eqnarray}
e^{(r)}_{n}e^{(s)}=\delta_{n,0}\ell \<e^{(r)},e^{(s)}\>{\bf 1}
=2\ell \delta_{n,0}(\delta_{r,s+1}+\delta_{r,s-1}){\bf 1}\ \ \ \mbox{ for }n\ge 0.
\end{eqnarray}
We also have
\begin{eqnarray}
e^{(r)}_{-1}e^{(r)}=0\ \ \ \mbox{ for }r\in \Z.
\end{eqnarray}

\bl{lsimplicity} Assume that $\ell$ is a nonzero complex number. Then
$V_{\hat{E}}(\ell,0)$ viewed as a $\Cl(L(E),\ell)$-module is
irreducible and $V_{\hat{E}}(\ell,0)$ viewed as a vertex
super-algebra is simple. \el

\begin{proof}  Let $k$ be any positive
integer. Set $E[k]={\rm span}\{e^{(r)}\;|\; -k\le r\le k+1\}$ (even
dimensional). One can show that the bilinear form of $E$, restricted to
$E[k]$, is non-degenerate. Let $V[k]$ denote the vertex super subalgebra of
$V_{\hat{E}}(\ell,0)$, generated by $E[k]$.  {}From \cite{ffr}, $V[k]$
is an irreducible $\Cl(L(E[k]),\ell)$-module. It is clear that $V[k]$ with $k\ge 1$ form an increasing filtration of
$V_{\hat{E}}(\ell,0)$. Then it follows that
$V_{\hat{E}}(\ell,0)$ is an irreducible $\Cl(L(E),\ell)$-module.
This particularly implies that $V_{\hat{E}}(\ell,0)$ is a simple vertex
super-algebra.
\end{proof}

For the rest of this section we consider $\ell=1$.
The following is straightforward:

\bl{laction-aut} For $m\in \Z$, there exists an automorphism
$\sigma_{m}$ of $V_{\hat{E}}(1,0)$, which is uniquely determined by
$$\sigma_{m}(e^{(r)})=e^{(m+r)}\ \ \ \mbox{ for }r\in \Z.$$
Furthermore, this gives rise to a group action of $\Z$ on $V_{\hat{E}}(1,0)$
by automorphisms. \el

Define a linear character of $\Z$:
\begin{eqnarray}
\chi_{p}: \Z\rightarrow \C^{\times}\ \mbox{ with } \chi_{p}(n)=p^{n} \mbox{ for }n\in \Z.
\end{eqnarray}
Now, we are in a position to present our first main result of this section.

\bt{tmain-p-1} Let $W$ be any restricted ${\mathcal{V}}ir_{p,-1}$-module.
Then there exists a $(\Z,\chi_{p})$-covariant $\phi$-coordinated quasi
$V_{\hat{E}}(1,0)$-module structure $Y_{W}(\cdot,x)$ which is
uniquely determined by the condition
 $$Y_{W}(e^{(r)},x)=T(p^{r}x)\ \ \mbox{ for }r\in \Z.$$ \et

\begin{proof}  As $E$ generates $V_{\hat{E}}(1,0)$ as a vertex super-algebra,
the uniqueness is clear. We now prove the existence. Set
$$U_{W}=\{ T(p^{m}x)\;|\; m\in \Z\}\subset \E(W).$$
{}From (\ref{erelation-generating-q-5}) we have
\begin{eqnarray}\label{eneed-proof}
&&\left(\frac{x_{1}+p^{s-r}x_{2}}{x_{1}-p^{s-r}x_{2}}\right)T(p^{r}x_{1})T(p^{s}x_{2})
-\left(\frac{x_{2}+p^{r-s}x_{1}}{x_{2}-p^{r-s}x_{1}}\right)T(p^{s}x_{2})T(p^{r}x_{1})\nonumber\\
&&\ \ \ \hspace{1cm} =-2\left(\frac{1+p}{1-p}\right)
\left(\delta\left(\frac{p^{s+1-r}x_{2}}{x_{1}}\right)
-\delta\left(\frac{p^{s-1-r}x_{2}}{x_{1}}\right)\right)
\end{eqnarray}
for $r,s\in \Z$. By multiplying both sides of (\ref{eneed-proof}) by a suitable polynomial we
get
\begin{eqnarray}\label{e5.23}
&&(x_{1}/x_{2}-p^{s+1-r})(x_{1}/x_{2}-p^{s-1-r})
(x_{1}/x_{2}+p^{s-r})T(p^{r}x_{1})T(p^{s}x_{2})\nonumber\\
&=&-(x_{1}/x_{2}-p^{s+1-r})(x_{1}/x_{2}-p^{s-1-r})
(x_{1}/x_{2}+p^{s-r})T(p^{s}x_{2})T(p^{r}x_{1}).\ \ \ \
\end{eqnarray}
{}From definition, $U_{W}$ is $\Gamma_{p}$-stable and from (\ref{e5.23}), $U_{W}$ is quasi $\S_{trig}$-local.
In view of Theorem \ref{tmain-phi}, $U_{W}$ generates a weak quantum vertex algebra $\<U_{W}\>_{e}$
with $W$ as a faithful $\phi$-coordinated quasi module.

We next show
that $\<U_{W}\>_{e}$ becomes a $\Cl(L(E),1)$-module by letting
$e^{(r)}(z)$ act as $Y_{\E}^{e}(T(p^{r}x),z)$ for $r\in \Z$.
To this end, we need to prove
\begin{eqnarray}\label{eneed-super}
&& Y_{\E}^{e}(T(p^{r}x),x_{1})Y_{\E}^{e}(T(p^{s}x),x_{2})+
Y_{\E}^{e}(T(p^{s}x),x_{2})Y_{\E}^{e}(T(p^{r}x),x_{1})\nonumber\\
&&\ \ \ \ =2(\delta_{r,s+1}+\delta_{s,r+1})x_{1}^{-1}
\delta\left(\frac{x_{2}}{x_{1}}\right)
\end{eqnarray}
for $r,s\in \Z$. With (\ref{e5.23}), by Proposition
5.3 of \cite{li-phi-module}, we have
\begin{eqnarray*}
P(e^{x_{1}-x_{2}})
Y_{\E}^{e}(T(p^{r}x),x_{1})Y_{\E}^{e}(T(p^{s}x),x_{2})
=-P(e^{x_{1}-x_{2}}) Y_{\E}^{e}(T(p^{s}x),x_{2})Y_{\E}^{e}(T(p^{r}x),x_{1}),
\end{eqnarray*}
where $P(x)=(x-p^{s+1-r})(x-p^{s-1-r})(x+p^{s-r})$.
As $p$ is not a root of $-1$, $e^{x_{1}-x_{2}}+p^{s-r}$ is
invertible in $\C[[x_{1},x_{2}]]$. By cancellation we get
\begin{eqnarray}\label{efirst-relation}
&&(e^{x_{1}-x_{2}}-p^{s+1-r})(e^{x_{1}-x_{2}}-p^{s-1-r})
Y_{\E}^{e}(T(p^{r}x),x_{1})Y_{\E}^{e}(T(p^{s}x),x_{2})\nonumber\\
&=&-(e^{x_{1}-x_{2}}-p^{s+1-r})(e^{x_{1}-x_{2}}-p^{s-1-r})
Y_{\E}^{e}(T(p^{s}x),x_{2})Y_{\E}^{e}(T(p^{r}x),x_{1}).\ \ \ \
\end{eqnarray}
If $r-s\ne \pm 1$,
$(e^{x_{1}-x_{2}}-p^{s+1-r})$ and $(e^{x_{1}-x_{2}}-p^{s-1-r})$ are also
invertible in $\C[[x_{1},x_{2}]]$, so that we furthermore obtain
\begin{eqnarray}
Y_{\E}^{e}(T(p^{r}x),x_{1})Y_{\E}^{e}(T(p^{s}x),x_{2}) =-
Y_{\E}^{e}(T(p^{s}x),x_{2})Y_{\E}^{e}(T(p^{r}x),x_{1}).
\end{eqnarray}
This proves that (\ref{eneed-super}) holds for $r\ne s\pm 1$.

Consider the case $r=s\pm 1$. Set $h(x)=\sum_{n\ge 1}\frac{1}{n!}x^{n-1}$.
Then $e^{x}-1=xh(x)$. If $r=s+1$, we have
$$(e^{x_{1}-x_{2}}-p^{s+1-r})(e^{x_{1}-x_{2}}-p^{s-1-r})
=(x_{1}-x_{2})h(x_{1}-x_{2})(e^{x_{1}-x_{2}}-p^{-2}),$$
where $h(x_{1}-x_{2})(e^{x_{1}-x_{2}}-p^{-2})$ is invertible in $\C[[x_{1},x_{2}]]$.
{}From (\ref{efirst-relation}) we get
\begin{eqnarray}
&&(x_{1}-x_{2})Y_{\E}^{e}(T(p^{r}x),x_{1})Y_{\E}^{e}(T(p^{s}x),x_{2})\nonumber\\
&=&-(x_{1}-x_{2})
Y_{\E}^{e}(T(p^{s}x),x_{2})Y_{\E}^{e}(T(p^{r}x),x_{1}).
\end{eqnarray}
If $r=s-1$, similarly we obtain the same relation. As $\<U_{W}\>_{e}$ as a weak quantum vertex algebra
is generated by $T(p^{m}x)$ for $m\in \Z$,
it follows that $\<U\>_{e}$ is actually a vertex super-algebra.

Recall that $P(x)=(x-p^{s+1-r})(x-p^{s-1-r})(x+p^{s-r})$.
For $r=s+1$, we have $P(x)=(x-1)(x-p^{-2})(x+p^{-1})$
and $P'(1)=(1-p^{-2})(1+p^{-1})$. With (\ref{e5.23}), by Lemma \ref{labuv} we have
$$T(p^{r}x)_{n}^{e}T(p^{s}x)=0\ \ \mbox{ for }n\ge 1$$
and
\begin{eqnarray*}
&&P'(1)T(p^{r}x)_{0}^{e}T(p^{s}x)\\
&=&\Res_{x_{1}}\left(\frac{1}{x_{1}-x}P(x_{1}/x)T(p^{r}x_{1})T(p^{s}x)+
\frac{1}{-x+x_{1}}P(x_{1}/x)T(p^{s}x)T(p^{r}x_{1})\right)\\
 &=&\Res_{x_{1}}x^{-1}(x_{1}/x-p^{-2})(x_{1}/x-p^{-1})\cdot \\
&&\ \ \ \ \cdot
\left(\left(\frac{x_{1}+p^{-1}x}{x_{1}-p^{-1}x}\right)T(p^{r}x_{1})T(p^{s}x)
-\left(\frac{p^{-1}x+x_{1}}{p^{-1}x-x_{1}}\right)T(p^{s}x)T(p^{r}x_{1})\right)\nonumber\\
&=&\Res_{x_{1}}x^{-1}(x_{1}/x-p^{-2})(x_{1}/x-p^{-1})\cdot \\
&&\ \ \ \ \cdot (-2)\left(\frac{1+p}{1-p}\right)
\left(\delta\left(\frac{x}{x_{1}}\right)
-\delta\left(\frac{p^{-2}x}{x_{1}}\right)\right)\\
&=&(1-p^{-2})(1-p^{-1})(-2)\left(\frac{1+p}{1-p}\right)\\
&=&2(1-p^{-2})(1+p^{-1})(\delta_{r,s+1}+\delta_{s,r+1}),
\end{eqnarray*}
which implies
$$T(p^{r}x)_{0}^{e}T(p^{s}x)=2\cdot 1_{W}.$$
It follows from the super-commutator formula that (\ref{eneed-super}) holds for $r=s+1$.
Similarly, one can show that (\ref{eneed-super}) holds for $r=s-1$.
Thus $\<U_{W}\>_{e}$ is a
$\Cl(L(E),1)$-module with $e^{(r)}(z)$ acting as
$Y_{\E}^{e}(T(p^{r}x),z)$ for $r\in \Z$. Furthermore, we have
$$e^{(r)}_{n}1_{W}=T(p^{r}x)_{n}^{e}1_{W}=0\ \ \mbox{ for }r\in \Z,\; n\in \N.$$
It follows from the construction of $V_{\hat{E}}(1,0)$ that there
exists a $\Cl(L(E),1)$-module homomorphism $\psi$ from
$V_{\hat{E}}(1,0)$ to $\<U_{W}\>_{e}$, sending ${\bf 1}$ to $1_{W}$.
That is,
$$\psi(Y(e^{(r)},z)v)=Y_{\E}^{e}(T(p^{r}x),z)\psi(v)\ \ \ \mbox{ for }r\in \Z,\ v\in V_{\hat{E}}(1,0).$$
Since $E$ generates $V_{\hat{E}}(1,0)$ as a vertex super-algebra,
it follows that $\psi$ is a homomorphism of
vertex super-algebras. As $W$ is a $\phi$-coordinated quasi module
for $\<U\>_{e}$,  $W$ becomes a $\phi$-coordinated quasi
 $V_{\hat{E}}(1,0)$-module through homomorphism $\psi$. Furthermore,
we have
$$Y_{W}(\sigma_{n}(e^{(r)}),x)=Y_{W}(e^{(r+n)},x)=T(p^{r+n}x)=Y_{W}(e^{(r)},p^{n}x)=Y_{W}(e^{(r)},\chi_{p}(n)x)$$
for $n,r\in \Z$. As $E$ generates $V_{\hat{E}}(1,0)$, by Lemma \ref{l-phicoordinated-module-induction}, it
is also $(\Z,\chi_{p})$-covariant. Therefore, $W$ is a $(\Z,\chi_{p})$-covariant $\phi$-coordinated quasi $V_{\hat{E}}(1,0)$-module.
\end{proof}

On the other hand, we have:

\bt{tmain-p-2} Let $(W,Y_{W})$ be a $(\Z,\chi_{p})$-covariant
$\phi$-coordinated quasi $V_{\hat{E}}(1,0)$-module. Then $W$ is a restricted
module for ${\mathcal{V}}ir_{p,-1}$ with
$T(x)=Y_{W}(e^{(1)},x)$. \et

\begin{proof} As $Y_{W}(e^{(1)},x)\in \E(W)$ from definition, we must prove
\begin{eqnarray}\label{eneed-toshow}
&&\left(\frac{x_{1}+x_{2}}{x_{1}-x_{2}}\right)
Y_{W}(e^{(1)},x_{1})Y_{W}(e^{(1)},x_{2})
-\left(\frac{x_{2}+x_{1}}{x_{2}-x_{1}}\right)Y_{W}(e^{(1)},x_{2})Y_{W}(e^{(1)},x_{1})
\nonumber\\
&=&-2\left(\frac{1+p}{1-p}\right)
\left(\delta\left(\frac{px_{2}}{x_{1}}\right)
-\delta\left(\frac{p^{-1}x_{2}}{x_{1}}\right)\right).
\end{eqnarray}
Note that for $n\in \Z,\ i\ge 0$, we have
$$(\sigma_{n}e^{(1)})_{i}e^{(1)}=e^{(n+1)}_{i}e^{(1)}=\delta_{i,0}\<e^{(n+1)},e^{(1)}\>{\bf 1}
=2\delta_{i,0}(\delta_{n,1}+\delta_{n,-1}){\bf 1}.$$
In view of Theorem \ref{tcommutator-quasi-phi-module}, we have
\begin{eqnarray}\label{ecomparison-1}
&&Y_{W}(e^{(1)},x_{1})Y_{W}(e^{(1)},x_{2})+Y_{W}(e^{(1)},x_{2})Y_{W}(e^{(1)},x_{1})\nonumber\\
&=&2\left(\delta\left(\frac{px_{2}}{x_{1}}\right)+\delta\left(\frac{p^{-1}x_{2}}{x_{1}}\right)\right),
\end{eqnarray}
noticing that $\chi_{p}$ is one-to-one.
{}From this we obtain
\begin{eqnarray*}
&&(x_{1}-px_{2})(px_{1}-x_{2})Y_{W}(e^{(1)},x_{1})Y_{W}(e^{(1)},x_{2})\\
&=&-(x_{1}-px_{2})(px_{1}-x_{2})Y_{W}(e^{(1)},x_{2})Y_{W}(e^{(1)},x_{1}),
\end{eqnarray*}
which implies
\begin{eqnarray}
(x_{1}-px_{2})(px_{1}-x_{2})Y_{W}(e^{(1)},x_{1})Y_{W}(e^{(1)},x_{2})\in \Hom(W,W((x_{1},x_{2}))).
\end{eqnarray}

Note that $e^{(1)}_{j}e^{(1)}=0$ for $j\ge -1$. Using the fact that $Y_{W}$ is a homomorphism
(from Proposition \ref{pquasi-representation}), we get
$$Y_{W}(e^{(1)},x)_{j}^{e}Y_{W}(e^{(1)},x)=Y_{W}(e^{(1)}_{j}e^{(1)},x)=0\ \ \mbox{ for }j\ge -1.$$
{}From Lemma \ref{ltruncation} we have
$$(x_{1}-x_{2})^{-1}(x_{1}-px_{2})(px_{1}-x_{2})Y_{W}(e^{(1)},x_{1})Y_{W}(e^{(1)},x_{2})
\in \Hom(W,W((x_{1},x_{2}))).$$
Thus
\begin{eqnarray}
(x_{1}-px_{2})(px_{1}-x_{2})Y_{W}(e^{(1)},x_{1})Y_{W}(e^{(1)},x_{2})&=&(x_{1}-x_{2})A(x_{1},x_{2}),\nonumber\\
(x_{1}-px_{2})(px_{1}-x_{2})Y_{W}(e^{(1)},x_{2})Y_{W}(e^{(1)},x_{1})&=&(x_{2}-x_{1})A(x_{1},x_{2})
\end{eqnarray}
for some $A(x_{1},x_{2})\in \Hom (W,W((x_{1},x_{2})))$.
Then
\begin{eqnarray}
&&Y_{W}(e^{(1)},x_{1})Y_{W}(e^{(1)},x_{2})+Y_{W}(e^{(1)},x_{2})Y_{W}(e^{(1)},x_{1})\nonumber\\
&=&\left(\frac{x_{1}-x_{2}}{(x_{1}-px_{2})(px_{1}-x_{2})}+\frac{x_{2}-x_{1}}{(x_{2}-px_{1})(px_{2}-x_{1})}\right)
A(x_{1},x_{2})
\nonumber\\
&=&\frac{1}{p+1}\left(x_{1}^{-1}\delta\left(\frac{px_{2}}{x_{1}}\right)
+x_{2}^{-1}\delta\left(\frac{px_{1}}{x_{2}}\right)\right)A(x_{1},x_{2}).
\end{eqnarray}
Comparing this with (\ref{ecomparison-1}) we get
\begin{eqnarray*}
&&x_{1}^{-1}\delta\left(\frac{px_{2}}{x_{1}}\right)A(x_{1},x_{2})=2(p+1)\delta\left(\frac{px_{2}}{x_{1}}\right),\\
&&x_{2}^{-1}\delta\left(\frac{px_{1}}{x_{2}}\right)A(x_{1},x_{2})=2(p+1)\delta\left(\frac{px_{1}}{x_{2}}\right).
\end{eqnarray*}
Using these identities we obtain
\begin{eqnarray*}
&&\left(\frac{x_{1}+x_{2}}{x_{1}-x_{2}}\right)
Y_{W}(e^{(1)},x_{1})Y_{W}(e^{(1)},x_{2})
-\left(\frac{x_{2}+x_{1}}{x_{2}-x_{1}}\right)Y_{W}(e^{(1)},x_{2})Y_{W}(e^{(1)},x_{1})
\nonumber\\
&=&\left( \frac{x_{1}+x_{2}}{(x_{1}-px_{2})(px_{1}-x_{2})}-\frac{x_{1}+x_{2}}{(px_{2}-x_{1})(x_{2}-px_{1})}\right)A(x_{1},x_{2})\\
&=&\frac{1}{p-1}\left( \frac{1}{x_{1}-px_{2}}-\frac{1}{px_{1}-x_{2}}+\frac{1}{px_{2}-x_{1}}-\frac{1}{x_{2}-px_{1}}\right)A(x_{1},x_{2})\\
&=&\frac{1}{p-1}\left( x_{1}^{-1}\delta\left(\frac{px_{2}}{x_{1}}\right)
-x_{2}^{-1}\delta\left(\frac{px_{1}}{x_{2}}\right)\right)A(x_{1},x_{2})\\
&=&-2\left(\frac{1+p}{1-p}\right)
\left(\delta\left(\frac{px_{2}}{x_{1}}\right)
-\delta\left(\frac{p^{-1}x_{2}}{x_{1}}\right)\right).
\end{eqnarray*}
Therefore, $W$ is a restricted
module for ${\mathcal{V}}ir_{p,-1}$ with
$T(x)=Y_{W}(e^{(1)},x)$.
\end{proof}

\end{document}